\documentclass[12pt]{amsart}
\usepackage{amsmath,latexsym,amssymb}
\begin{document}

\newtheorem{definition}{Definition}[section]
\newtheorem{theorem}{Theorem}[section]
\newtheorem{proposition}[theorem]{Proposition}
\newtheorem{lemma}[theorem]{Lemma}
\newtheorem{remark}{Remark}
\newtheorem{corollary}[theorem]{Corollary}
\newtheorem{question}{Question}
\newtheorem{example}{Example}
\newtheorem{notation}{Notation}[section]
\newtheorem{claim}{Claim}[section]
\def\cl{\begin{claim}\upshape}
\def\ecl{\end{claim}}
\def\prcl{\par\noindent{\em Proof of Claim: }}
\def\bem{\begin{remark}\upshape}
\def\ebem{\end{remark}}
\def\nota{\begin{notation}\upshape}
\def\enota{\end{notation}}
\def\defn{\begin{definition}\upshape}
\def\edefn{\end{definition}}
\def\thm{\begin{theorem}}
\def\ethm{\end{theorem}}
\def\lmm{\begin{lemma}}
\def\elmm{\end{lemma}}
\def\qed{\hfill$\quad\Box$}
\def\pr{\par\noindent{\em Proof: }}
\def\prbeh{\par\noindent{\em Proof of Claim: }}
\def\kor{\begin{corollary}}
\def\ekor{\end{corollary}}
\def\frag{\begin{question}\upshape}
\def\efrag{\end{question}}
\def\prop{\begin{proposition}}
\def\eprop{\end{proposition}}
\def\beh{\par\noindent{\bf Claim: }\em}
\def\ebeh{\par\noindent\upshape}
\def\bsp{\begin{example}}
\def\ebsp{\end{example}}
\def\a{{\sigma}}
\def\si{{\sigma}}
\def\F{{ \mathbb F}}
\def\N{{\mathbb N}}
\def\Z{{\mathbb Z}}
\def\Q{{\mathbb Q}}
\def\L{{\mathcal L}}
\def\P{{\mathcal P}}
\def\I{{\mathcal I}}
\def\B{{\mathcal B}}
\def\C{{\mathcal C}}
\def\R{{\mathcal R}}
\def\IR{{\mathbb R}}
\def\M{{\mathcal M}}
\def\I{{\mathcal I}}
\def\O{{\mathcal O}}
\def\Ma{{\mathfrak m}}
\def\U{{\Upsilon}}
\title{Separably closed fields and contractive Ore modules}
\date{\today}
\author{Luc B\'elair}
\address{\hskip-\parindent
Luc B\'elair\\
D\'epartment de math\'ematiques\\
Universit\'e du Qu\'ebec-UQAM, C.P. 8888 succ.
centre-ville \\ Montr\'eal, Qu\'ebec, H3C 3P8.}
\email {belair.luc@uqam.ca}
\author{Fran\c coise Point$^1$}
\address{\hskip-\parindent
Fran\c coise Point\\
Institut de Math\'ematique\\
Universit\'e de Mons, Le Pentagone\\
20, place du Parc, B-7000 Mons, Belgium}
\email {point@math.univ-paris-diderot.fr}
\thanks{\\
${}^1$Research Director at the "Fonds de la Recherche
Scientifique FNRS-FRS", partially funded by the FNRS-FRS Research Credit 14729555.}

\begin{abstract}
We consider valued fields with a distinguished contractive map as valued modules over the Ore ring of difference operators. We prove quantifier elimination for separably closed valued fields with the Frobenius map, in the pure module language augmented with functions yielding components for a $p$-basis and a chain of subgroups indexed by the valuation group.
\end{abstract}
\maketitle
\section{Introduction}

\par Let $K$ be a valued field of prime characteristic $p$, and let $Frob_p$ denote the Frobenius map $x\mapsto x^p$, and $v$ the valuation map. In \cite{R}, Rohwer studied the additive structure $(K,+,Frob_p)$ in a formalism taking into account the valuation through the chain of subgroups $V_\delta=\{ x : v(x)\ge \delta\}$, and he proved model-completeness for such models as  $K={\mathbb{F}}_p((T))$ and $K=\widetilde{\mathbb{F}}_p((T))$, $\widetilde{\mathbb{F}}_p$ being the algebraic closure of ${\mathbb{F}}_p$. We recall that the corresponding full theory of valued fields has been studied (see \cite{Az}), but is in general very far from being fully understood (see e.g. \cite{DS}), in particular for the above two examples. In \cite{BP},  we investigated the additive theory of valued fields but with a distinguished isometry (at the opposite of the Frobenius map) and we could obtain results similar to Rohwer's, even at the level of quantifier elimination for such models as $K=\widetilde{\mathbb{F}}_p((T))$ with the isometry $\sigma(\sum a_iT^i)=\sum a_i^pT^i$. In contrast with Rohwer, our starting point does not address directly the structure of some specific classes of definable sets, but is in the spirit of classical elimination of quantifiers algorithms in the theory of modules. In this paper, we show that our methods can be applied to the Frobenius map for separably closed valued fields (Proposition \ref{EQsep}), a case not covered by Rohwer (see below, \S 4).  In order to describe the theory of modules over the Ore ring of difference operators, we will use the formalism of $\lambda$-functions introduced by G. Srour (\cite{Sr}, see also \cite{D}, \cite{H}), and follow the approach undertaken, for instance, in \cite{DDP1}, \cite{DDP2} and \cite{P}. Finally let us mention that new results have been obtained by G. Onay on the model theory of valued modules (\cite{On}), both in the isometric case and the contractive case (the Frobenius map case). 

\par We mostly use the notation of \cite{BP}, with some slight modifications.

\section{Rings of power series as modules}
\par 
Let $D$ be a ring with a distinguished endomorphism $\sigma$ and let $A_{0}:=D[t;\sigma]$ the corresponding skew polynomial ring with the commutation rule $a.t=t.a^{\sigma}$ (\cite{C} chapter 2). Recall that any element of $A_{0}$ can be written uniquely as $\sum_{i=0}^n t^i.a_{i}$, with $a_{i}\in D$ (\cite{C} Proposition 2.1.1 (i)) and one has a degree function $deg : A_0\setminus\{0\} \to \mathbb N$ sending  $\sum_{i=0}^n t^i.a_{i}$ with $a_{n}\neq 0$, to $n\in \N$, and  with the convention that $deg(0)=-\infty<\N$. We can consider $D$ as an $ A_{0}$-module, by interpreting scalar multiplication by $t$ by the action of $\sigma$ on $D$.
\par In addition, we will assume that $D$ is a right Ore domain, namely that for all nonzero $a,b$ there are nonzero $c, d$ such that $\;a.c=b.d$, and that $\sigma$ is injective, which yields that $A_{0}$ has no zero divisor (\cite{C} Proposition 2.1.1 (ii)). 
\par Under these assumptions, $ A_{0}$ satisfies 
the {\it generalized right division algorithm}:
for any $q_1(t), \;q_2(t)\in  A_{0}$ with $deg(q_1)\geq deg(q_2)$, there exist $a\in D-\{0\},\;d\in \N$ and
$c$, $r$ in $ A_{0}$ with $deg(r)<deg(q_2)$ such that $q_1.a^d=q_2.c+r$ (see e.g. \cite{BP}, Lemma 2.2). 
\par Since $D$ is a right Ore domain, it has a right field of fractions $K$ and we denote the extension of $\sigma$ on $K$ by the same letter. Let $A:=K[t;\sigma]$.
\par Let $K^{\sigma}$ the subfield of $K$ consisting of the image of $K$ under $\sigma$. We will consider $K$ as a $K^{\sigma}$ vector-space and we will fix a basis $\C$. We will call such basis a {\it  $\sigma$-basis}.
\par Moreover we will assume that $\C$ can be chosen in $D$ and that
any element of $D$ has a decomposition along that basis with coefficients in $D^{\sigma}$. This is the case for instance if $K$ has characteristic $p$ and $\sigma$ is the Frobenius endomorphism.
\par For simplicity, we will assume that $\C$ is finite, that it contains $1$ and we present $\C$ as a finite tuple of distinct elements $(1=c_0,\cdots,c_{n-1})$.  However the infinite-dimensional case is not essentially different (see \cite{P}).
\par Later we will need both $A_{0}$ and $A$, but for the moment we will denote by the letter $A$ a skew polynomial ring of the form $D[t;\sigma]$, where $D$ satisfies the above hypothesis (which encompasses the case where $D$ is a field).
\par We will consider $A$-modules $M$ which have a direct sum decomposition as follows : $M=\oplus_{i=0}^{n-1} M.tc_{i}$. We will add new unary function symbols $\lambda _{i}$, $i\in n=\{0,\cdots,n-1\}$ to the usual language of $A$-modules in order to ensure the existence of such decomposition in the class of $A$-modules we will consider. These functions will be additive and so we will stay in the setting of abelian structures (see for instance \cite{Pr}).


\defn\label{def3} Let $L_{A}:=\{+,-,0,\cdot a; a\in A\}$ be the usual language of $A$-modules. Let $\lambda _{i}$, $i\in n=\{0,\cdots,n-1\}$, be new unary function symbols. Let 
$\L_A =L_A\cup
\{\lambda_{i}; i\in n\}$, 
and let $T_{\sigma}$ be the following $\L_A$-theory:
\begin{enumerate}
\item  
 the $L_{A}$-theory of all right $A$-modules
\item 
$\forall x\;\;(x=\sum_{i\in n}\lambda _{i}(x)\cdot t c_i)$
\item $\forall x \forall (x_{i})_{i\in n}(x=
\sum_{i\in n} x_{i}\cdot t c_i\,\rightarrow 
\bigwedge_{i\in n} x_{i}=
\lambda_i(x)$).
\end{enumerate}
\edefn
\par Note that $D$ is a model ot $T_\sigma$, when viewed as an $A$-module as before.
\par We will need later that the functions $\lambda_i$, $i\in n$ are defined in any model of $T_{\sigma}$ by the following $\L_{A}$-formula:
$\lambda_i(x)=y$ iff  
($\exists y_0\cdots \exists y_{n-1} \;\; x=\sum_{j\in n} y_j\cdot t c_j$ and
 $ y_i=y$) iff ($\forall y_0\cdots \forall y_{n-1} \;\; x=\sum_{j\in n} y_j\cdot t c_j\;\;\rightarrow\;\;y_i=y$).
\par Such theories $T_{\sigma}$ have been investigated in \cite{DDP1},  \cite{DDP2} and \cite{P}, when $D$ is a separably closed field of characteristic $p$ and $\sigma$ is the Frobenius endomorphism. Let us recall some of the terminology developed there.
\nota An element $q(t)$ of $A$ is {\it $\sigma$-separable}
if  $q(0)\neq 0$. 
In writing down an element of $A$, we will allow ourselves
to either write it as $q$ or $q(t)$ when stressing the fact that it is a 
polynomial in $t$.
\enota
\par In order to reduce divisibility questions to divisibility by separable polynomials, it is convenient to introduce the following notation.
\nota\label{lambda} (See Notation 3.2, Remark 2 and section 4 in \cite{DDP1}.)
\par Given $q\in A$ , we will define $\sqrt[\sigma]{q}$ and 
$q^{\sigma}$. First, for 
$a=\sum_{i} a_{i}^{\sigma}c_{i}\in D$, where the elements $a_{i}$
belong to $D$ and $c_i$'s to $\C$, 
set $a^{1/{\sigma}}:=\sum_{i} a_{i}c_{i}$.
Observe that $(a^{\sigma})^{1/{\sigma}}=a$, but unless 
$a\in A^{\sigma}$, $(a^{1/{\sigma}})^{\sigma}$ and $a$ are distinct. 
Then, for $q=\sum_{j=0}^{n}t^{j}a_{j}\in A$ with $a_{j}\in D$, 
set $\sqrt[\sigma]{q}:=\sum_{j=0}^{n}t^{j}a_{j}^{1/{\sigma}}$.
We also define $q^{\sigma}$ as $\sum_{i=0}^{n}t^{j}a_{j}^{\sigma}$, we have $tq^{\sigma}=qt$. 
\par Iteration of $\sqrt[\sigma]{\ }$ $m$ times is denoted by $\sqrt[\sigma^{m}]{\ }$. Let 
$a=\sum_{i=0}^{n-1} a_{i}^{\sigma}c_{i}\in D$, where $a_{i}\in D$ and $c_i\in\C$. Decompose each $a_{i}$ along the basis $\C$, $a_{i}=\sum_{j=0}^{n-1} a_{ij}^{\sigma}c_{j},$ so $a_{i}^{\sigma}=\sum_{j=0}^{n-1} a_{ij}^{\sigma^2}c_{j}^{\sigma}$ and $a=\sum_{i,\;j=0}^{n-1} a_{ij}^{\sigma^2}c_{j}^{\sigma}c_{i}$. More generally, $a=\sum_{\bar d\in n^m} a_{\bar d}^{\sigma^m}c_{\bar d}$, where $\bar d:=(d_{1},\cdots,d_{m})\in n^m$, $c_{\bar d}:=c_{d_{1}}^{\sigma^{m-1}}\cdots c_{d_{m}}$. 
\par Given $q\in A$, we write it as $q=\sum_i q_i c_i$ with $q_i=\sum_j t^j a_{ij}^{\sigma}$, $a_{ij} \in D$. 
Therefore, we have that $\sqrt[\sigma]{q_i}=\sum_j t^j a_{ij}$, so
$$t q=\sum_i \sqrt[\sigma]{q_i} t c_i.$$
Indeed, $\sum_i \sqrt[\sigma]{q_i} t c_i=\sum_i \sum_j t^j a_{ij} t c_i=
\sum_i \sum_j t^{j+1} a_{ij}^{\sigma} c_i=t q.$
\par Similarly, $t^m q=\sum_{\bar d\in n^m}\sqrt[\sigma^{m}] q_{\bar d}t^m c_{\bar d}$.
\enota

\par For example, let $F$ be a field of characteristic $p$, $D=F[[x]]$ and $\sigma$ be the Frobenius map on $D$. The notion of $\sigma$-separable polynomials coincides with the notion of $p$-polynomials $\sum_{i=0}^m a_{i} x^{p^{m-i}}$, with $a_{m}\neq 0$, introduced by O. Ore \cite{Ore2} (see also \cite{DDP1}). In case $F$ is perfect, a $\sigma$-basis for $F[[x]]$ is $\{1,x,\cdots,x^{p-1}\}$.
In general, let $\B$ be a $p$--basis of $F$. Then 
$D^{\sigma}$ is equal to $F^{\sigma}[[x^p]]$ and there is a direct sum
decomposition of $D$ as 
$\oplus_{c_i\in\{\B,x.\B,\cdots,x^{p-1}.\B \}}  D^{\sigma}\cdot c_i$. 


\medskip



\par We now assume that $(K,v)$ is a valued field with valuation ring $\O_{K}$, maximal ideal $\Ma_{K}$ and residue field $\bar K$. Let $K^{\times}=K\setminus\{0\}$, we denote the value group $v( K^{\times})$ by $\Gamma$.
\par We will set $\bar a= a+\Ma_{K}$, the image of $a$ under the residue map from $\O_{K}$ to $\bar K$. Moreover, $K$ is endowed with an endomorphism $\sigma$ which is (valuation) increasing on $\O_{K}$ 
and strictly increasing on $\Ma_{K}$, namely, the following holds for all $a\in K$ : $$(v(a)\geq 0\rightarrow (v(a^{\sigma})\geq v(a))\;\&\;(v(a)>0\rightarrow (v(a^{\sigma})> v(a))).$$ 
This implies that $\sigma$ is an isometry on the elements of valuation zero and strictly decreasing on the elements of negative valuation. 
 In particular $\sigma$ induces an endomorphism $\sigma_{v}$ of $(\Gamma,+,\leq,0)$  defined by $$\sigma_{v}(v(a)):=v(\sigma(a)).$$ Note that $\sigma_{v}$ is injective. 
 We will denote the image by $\sigma_{v}$ of an element $\gamma\in \Gamma$  either by $\gamma^{\sigma_{v}}$ or by $\sigma_{v}(\gamma)$. In the example above where $\sigma$ is the Frobenius map, we have $\sigma_v(\gamma)=p\gamma.$ This action induced by $\sigma$ on the value group makes it a {\it multiplicative ordered difference abelian group} in the terminology used by K. Pal (\cite{Pa}), who investigated the model theory of such structures arising in the context of valued difference fields.
\par Let $A:= K[t;\sigma]$, $ A_{0}=\O_{K}[t; \sigma]$. We extend the residue map to $A_0$ by sending $q(t)=\sum_{j} t^j a_{j}$ to $\bar q(t):=\sum_{j} t^j \bar a_{j}$. We denote by $\I$ the set of elements of $A_{0}$ which have at least one coefficient of valuation $0$ (or equivalently $\bar q(t)\ne \bar 0$). Note that unlike the case where $\sigma$ is an isometry of $K$, one cannot extend the valuation $v$ of $ K^{\times}$ to $A^{\times}$ or to $A_{0}^{\times}$, but the product of two elements of $\I$ still belongs to $\I$.
\medskip
\nota 
For $q(t)\in A$ and $\mu\in K-\{0\}$, denote by
$q^{\mu}(t)$ the element of $A$ equal to $\mu.q(t).\mu^{-1}$. So if $q(t)=\sum_{i} t^i a_{i}$, $a_{i}\in K$, $q^\mu(t)=
\sum_{i} t^i \mu^{\sigma^i}a_{i}\mu^{-1}.$
\enota 
\par Note that if $q(t)\in  A_{0}$ and $\mu\in \O_{K}$, then $q^\mu(t)\in  A_{0}$.
\medskip

\section{Valued modules}

We keep the same notation as in the previous section with a fixed $(K,v,\sigma)$, $\Gamma=v(K^\times)$ endowed with the induced endomorphism $\sigma_{v}$, $A$ the skew polynomial ring $K[t;\sigma]$ etc.

\par We will define the notion of $\sigma$-valued $A$-modules, or simply valued $A$-modules. Notions of valued modules occur in various places with many variations, see for instance \cite{Co} or \cite{L}(\S2). The following generalizes the notion in \cite{BP}.
\defn \label{Tw} {\rm (Cf. \cite{F}, \cite{Dries}, \cite{BP})}
 A {\em valued $A$-module} is a two-sorted structure 
$(M,(\Delta\cup\{+\infty\}, \le, 0_\Delta, +\gamma; \gamma\in \Gamma),w)$, where 
$M$ is an $A$-module, $(\Delta \cup \{+\infty\}, \le)$ is a totally ordered set for which $+\infty$ is a maximum, $0_\Delta\in \Delta$ a distinguished element, $+\gamma$ is an action of $\gamma\in \Gamma$ on $\Delta$, 
and $w$ is a map $w:M\to \Delta\cup\{ +\infty \}$ such that 

\begin{enumerate}
\item for all $\delta, \delta_1,\delta_2\in\Delta$, if $\delta_1\le\delta_2$ then $\delta_1+\gamma\le \delta_2+\gamma$, for each $\gamma\in \Gamma$, and $\delta+\gamma_1<\delta +\gamma_2$, for each $\gamma_1<\gamma_2\in \Gamma$; \label{actionG}\footnote{This axiom should replace axiom (1) in the definition given in \cite{BP}.}
\item for all $m_1, m_2\in M$, $w(m_1+m_2)\ge
\min\{w(m_1), w(m_2)\}$, and $w(m_1)=+\infty$ iff
$m_1=0;$
\item for all $m_{1}, m_{2}\in M$, $w(m_{1})<w(m_{2})$ iff $w(m_{1}\cdot t)<w(m_{2}\cdot t);$ \label{action_de_t}
\item for all $m\in M-\{0\}$, $w(m\cdot \mu)=w(m)+v(\mu)$, for each $\mu\in K^{\times}$.
\end{enumerate}
We denote the corresponding two-sorted language by $L_{w}$ and the corresponding theory by $T_{w}$. We will sometimes write $0$ instead of $0_\Delta$, for ease of notation.
\edefn
\par Taking $M=K$ and $w=v$, $\Delta=\Gamma$, then $K$ is a
valued
$A$-module with $t$ acting as $\sigma$ and $\Gamma$ acting on itself by translation.
\medskip
\par From the axioms above, we deduce as usual the 
following properties : $w(m)=w(-m)$, and if
$w(m_{1})<w(m_{2})$, then
$w(m_{1}+m_{2})=w(m_{1})$.
\medskip 
\par Note also, from axiom (1), that for each $m_{1},\;m_{2}\in M$ and $\mu\in K^{\times}$, if $w(m_1)\neq w(m_2)$ implies $w(m_1.\mu)\neq w(m_2.\mu)$.
\medskip
\par Note that $(w(M),\leq,0_{\Delta})$ is a substructure of $(\Delta\cup\{+\infty\},\le,0_\Delta)$, and that
$t$ induces an endomorphism $\tau$ on $(w(M),\leq,0_\Delta)$ defined by $w(m.t)=\tau(w(m))$. It is well-defined since by axiom (\ref{action_de_t}),  if $w(m_{1})=w(m_{2})$ then $w(m_{1}\cdot t)=w(m_{2}\cdot t)$.
\par From now on, we will impose a growth condition on the action of $t$ by introducing the following additional structure on $\Delta$, assuming now that $w$ is surjective. This will induce in particular that the action of $t$
on the corresponding class of valued $A$-modules will be {\it uniform}, with a compatibility condition between the action of $(\Gamma,\sigma_{v})$ and the action of $\tau$.

\defn \label{langageLV} Let $(\Delta,\le,0_{\Delta},\tau,+\gamma; \gamma\in \Gamma)$ be a totally ordered set with a distinguished element $0_\Delta$, $+\gamma$ an action of $\gamma\in\Gamma$  
 on $\Delta$, and $\tau$ a fixed endomorphism of $(\Delta, \le)$.
 \par We assume that the action of $\Gamma$ on $\Delta$ is transitive, and for all $\delta, \delta_1,\delta_2\in\Delta$, if $\delta_1\le\delta_2$ then $\delta_1+\gamma\le \delta_2+\gamma$, for each $\gamma\in \Gamma$, and $\delta+\gamma_1<\delta +\gamma_2$, for each $\gamma_1<\gamma_2\in \Gamma$.
\par The endomorphism $\tau$ satisfies the conditions, viz. : $\delta_1<\delta_2 \to \tau(\delta_1)<\tau(\delta_2), \tau(0_\Delta)=0_\Delta, \delta>0_\Delta \to \tau(\delta)>\delta, \delta<0_\Delta \to \tau(\delta)<\delta, $ and 
finally a compatibility condition between the action of $\sigma_{v}$ on $\Gamma$ and the action of $\tau$: for all $\gamma\in \Gamma$ we have $\tau(\delta+\gamma)=\tau(\delta)+\gamma^{\sigma_v}$.
\smallskip
\par Let us denote the corresponding language by $\L_{\Delta,\tau}$ and the corresponding theory by $T_{\Delta,\tau}$.
\edefn

Let $L_{w,\tau}:=L_w\cup \L_{\Delta,\tau}$; we will consider the
the class of $L_{w,\tau}$-structures \\$(M,(\Delta\cup\{+\infty\},\le,\tau,0_\Delta,+\gamma; \gamma\in \Gamma),w)$ satisfying the following properties  : 
\begin{enumerate}
\item[(1)] $(M,(\Delta\cup\{+\infty\},\le,0_\Delta,+\gamma; \gamma\in \Gamma),w)$ is a valued $A$-module;
\item[(2)] $(\Delta,\le,\tau,0_\Delta,+\gamma; \gamma\in \Gamma)$ is a model of $T_{\Delta,\tau}$;
\item[(3)] $w(m\cdot t)=\tau(w(m))$.
\end{enumerate}
\par Note that if $M\in \Sigma_{w,\tau}$, then if $w(m)>0_{\Delta}$, then $w(m\cdot t)>w(m)$, if $w(m)=0_{\Delta}$, then $w(m\cdot t)=w(m)$, and if $w(m)<0_{\Delta}$, then $w(m\cdot t)<w(m)$. Moreover, letting $\Gamma^+:=\{\gamma\in \Gamma:\;\gamma>0\}$, if $\gamma\in \Gamma^+$, then $\sigma_{v}(\gamma)>\gamma$. (Indeed, $0_{\Delta}+\gamma>0_{\Delta}$, so $\tau(0_{\Delta}+\gamma)>0_{\Delta}+\gamma$. By the compatibility condition, $\tau(0_{\Delta}+\gamma)=0_{\Delta}+\sigma_{v}(\gamma)$ and so $\sigma_{v}(\gamma)>\gamma$.)


\section{Abelian structures}



\par In order to stay into the setting of abelian
structures, we will use a less
expressive language. This language was used by T. Rohwer while considering the field of Laurent series over the prime field $\mathbb{F}_p$ with the usual Frobenius map $y\mapsto y^p$ (\cite{R}). Instead of the two-sorted structure $(M,\Delta,w)$, where $M$ is a valued $A$-module, he considered the one-sorted abelian structure
$(M, (M_{\delta})_{\delta\in \Delta}),$ where $M_\delta=\{ x\in M : w(x)\ge \delta\}$. Similarly, given a valued $A$-module $M$ with a $\sigma$-basis, we will add the functions $\lambda_i$ and consider the one-sorted abelian structure $(M, (\lambda_i)_{i\in n}, (M_\delta)_{\delta\in \Delta})$. We will put additional hypotheses on $\Delta$ and $\Gamma$ (the value group of $K$).

\par We will consider theories of abelian structures satisfying strong divisibility properties. The basic example is the separable closure of $\tilde \F_{p}((T))$. 
Note that this example is not covered by Rohwer, as we will indicate below, following Corollary \ref{actiont}. 

\defn Let $(\Delta,\le,\tau,+,0_\Delta,+\infty)$ be a model of $T_{\Delta,\tau}$. We set the language $\L_V:=\L_{A}\cup\{V_{\delta} : \delta\in \Delta\}$, where $V_{\delta}$ is a unary predicate. 
\par Let $T_{V}$ be the $\L_V$-theory with the following axioms:
\begin{enumerate}
\item $T_{\sigma}$;
\item $\forall m \; (V_{\delta_1}(m) \rightarrow V_{\delta_2}(m)),$ whenever $\delta_1\le\delta_2$;
\item $\forall m_{1}\;\forall m_{2}\;\;
(V_{\delta}(m_{1})\;\&\;V_{\delta}(m_{2}) 
\rightarrow V_{\delta}(m_{1}+m_{2}))$;
\item $\forall m\;\;(V_{\delta}(m)\rightarrow V_{\delta+v(\mu)}(m.\mu)),$ 
where $\mu\in  K^{\times}$;
\item $\forall m\;(V_{\delta}(m)\leftrightarrow V_{\tau(\delta)}(m.t))$.
\end{enumerate}
\edefn
\par If $(M,\Delta,w)$ is a valued
$A$-module with a $\sigma$-basis, and we let $\M_{V}=(M,+,0, (.r)_{r \in A},$ $(\lambda_{i})_{i\in n}, (M_{\delta})_{
\delta\in \Delta})$, then $\M_{V}$ is a model of $T_V$, where $V_\delta$ is interpreted as $M_\delta$.
\par The structure $\M_{V}$ is an abelian structure and 
one gets as in the classical case of (pure) modules
that any formula is equivalent to a boolean
combination of positive primitive formulas (p.p.) and index sentences
(namely, sentences telling the index of two 
p.p.-definable subgroups of the domain of $\M_{V}$ in one another 
(see \cite{P})) and this p.p. elimination is uniform in
the class of such structures. 
\par Note that the pure module theory of separably closed fields has quantifier elimination in the presence of the functions $\lambda_i$ (\cite{DDP1}).

We want to axiomatize a class of abelian structures which contains the class of valued separably closed fields. Note that the theory of valued separably closed fields has been shown to be model-complete in the language of valued fields augmented with predicates expressing $p$-independence (\cite{Del}) and to admit quantifier elimination in the language of valued fields augmented with the $\lambda_{i}$ functions (\cite{Ho}).
\par In the remainder of the section, we will recall certain properties of separably closed fields viewed as modules over the corresponding skew polynomial rings and we will formalise them in order to axiomatize the class of modules we will be working with. 

\prop\label{factorization} Assume $K$ is separably closed and let $q(t)\in I$. Then there exists a factorization of $q(t)$ into linear factors belonging to $\I$.
\eprop
\pr This follows from  \cite[chapter I, theorem 3]{Ore}. We give a proof here.
\par Let $q(t):=\sum_{i=0}^d t^{d-i} a_{i}$. First we factorize $q(t) a_{d}^{-1}=(t-f_{1}) \cdots (t-f_{d})$ with $f_{i}\in K$. 
 We apply the Euclidean algorithm in $K[t;\sigma]$ and so there exists $q_{1}(t)$ such that $q(t)=(t-f) q_{1}(t)+a$, for some $a\in K$ and $q_{1}(t)\in K[t;\sigma]$.
We want to show that we can choose $f$ such that $a=0$.
Write $q(t)=\sum_{i=0}^d t^{d-i} a_{i}$ and $q_{1}(t)=\sum_{j=0}^{d-1} t^{d-1-j} b_{j}$. Then we calculate $(t-f) q_{1}(t)$ and we express that $a=0$. We obtain that $f$ has to be a root of some separable polynomial with coefficients $a_{d},\cdots, a_{0}$. Namely, we get $a_{0}=b_{0}$, $a_{1}=b-1-f^{p^d} b_{0}=b_{1}-f^{p^d} a_{0}$, $a_{2}=b_{2}-f^{p^{d-2}}(a_{1}+f^{p^{d-1}} a_{0})$, and finally
$a_{d}=-f b_{d-1}=-f a_{d-1}-f^{p+1} a_{d-2}-\cdots-f^{1+p+p^2+\cdots+p^{d-1}} a_{0}$.
 \par So we have $q(t) a_{d}^{-1}=(t-f_{1}) \cdots (t-f_{d})$ with $f_{i}\in K$. If $f_{1}\notin \O_{K}$, then write
$(t-f_{1})=((t f_{1}^{-1})-1) f_{1}$.
Proceeding successively, we obtain a factorization of $q(t)$ into linear factors of the form $(t f_{i}'-1)$ or $(t-f_{i}')$ with $f_{i}'\in \O_{K}$, together with a constant factor of the form $f a_{d}, f\in K$.
\par So it remains to show that if $g(t)\in \I$, then both $g(t) (t f-1)$ and $g(t) (t-f)$ belong to $\I$, with $v(f)\geq 0$. W.l.o.g. we assume that in the first case $v(f)>0$.
\par Let $g(t)=\sum_{j=0}^n t^{n-j} b_{j}\in \I$ and let $\ell$ be minimal such that $v(b_{\ell})=0$. Calculate the coefficient of $t^{\ell}$ in $g(t) (t f-1)$. It is equal to $(-b_{\ell}+f b_{\ell-1})$ and of valuation $0.$ In the second case, we calculate the coefficient of $t^{\ell+1}$ in $g(t) (t-f)$, it is equal to $(b_{\ell}-f b_{\ell+1})$ and of valuation $0$.
\par \qed
\kor Suppose $K$ is separably closed and that given any $a, n_1,n_2\in\O_K, a\ne 0$, there exists $m_1, m_2\in\O_{K}$ such that $m_1\cdot (t-a)=n_1$, and $m_2\cdot (t a-1)=n_2.$ Then given any separable polynomial $q(t)\in \I$ and any element $n$ of $\O_{K}$ there exists $m\in\O_{K}$ such that $m\cdot q(t)=n$. Moreover, if $v(n)=0$, then there exists $m$ with $v(m)=0$ and $m\cdot q(t)=n$.
\ekor
\pr By the above proposition, it suffices to prove it for linear factors of the form $(t-a)$ or $(t a-1)$ with $a\in \O_{K}$. Let $n\in \O_{K}$. Then by assumption, there exists $m\in \O_K$ such that $n=m\cdot (t-a)$. 
If $v(n)=0$ and $v(m)>0$, then $\min\{v(m\cdot t), v(m\cdot a)\}>0$, a contradiction. Now let $m\in \O_K$ such that $m\cdot (t a-1)=n$. 
If $v(n)=0$, then necessarily $v(m)=0$.
\qed
\medskip
\lmm\label{angular} Assume that $\sigma_{v}$ induces multiplication by $p$ on $\Gamma$ and that $\Gamma$ is $p$-divisible.
Then, for any $\delta\in \Gamma$, and for any finite subset of $\Gamma$, $\{\gamma_{i}:\;0\leq i\leq d\}$, there exists $\mu\in \Gamma$ such that $\delta=\min_{0\le i\le d} \, p^i \mu+\gamma_{i}.$
\elmm
 \pr 
 Consider the lines $y=p^i x+\gamma_{i}$, $0\leq i\leq d$.
 They intersect the horizontal line $y=\delta$ in $d+1$ points (not necessarily distinct). Let $(\mu_{i},\delta)$ be the intersection points of the lines $y=p^i x+\gamma_{i}$ with the horizontal line $y=\delta$. Let $\mu$ be maximal such. Let $\delta_{i}=p^i \mu+\gamma_{i}$, $0\leq i\leq d$. So $(\mu,\delta_{i})$ are the intersection points of the other lines with the vertical line $x=\mu$. If for some $i$, $\delta_{i}<\delta$, then on the line $y=p^i x+\gamma_{i}$, we would have also the point
 $(\mu_{i},\delta)$ with $\mu_{i}<\mu$. But then the angular component of that line would be equal to $\frac{\delta_{i}-\delta}{\mu-\mu_{i}}$, which is negative, a contradiction. \qed
 
\defn We will say that $(\Gamma,+,-,0,\leq,\si_{v})$ is ordered linearly closed (o.l.-closed) 
if given any $\gamma_{0}, \ldots, \gamma_d, \gamma_{\delta}\in \Gamma$, there exists $\gamma_{\mu}\in \Gamma$ such that 
$\gamma_\sigma=\min_j \sigma^j_v(\gamma_\mu)+\gamma_j$. Note that such a $\gamma_{\mu}$ is unique since each $\sigma_{v}$ and each translation $+\gamma$ are strictly increasing on $\Gamma$.
\edefn

For example, if $K$ is separably closed of characteristic $p$ and $\sigma$ is the  Frobenius map, then $\sigma_v(\gamma)=p\gamma$, and $\Gamma=vK$ is divisible. So this ensures $\Gamma$ is o.l.-closed as we saw in Lemma \ref{angular}.
\par Because of the compatibility condition between the action of $\sigma_{v}$ on $\Gamma$ and the action of $\tau$ on $\Delta$, the fact that  $\Gamma$ is o.l.-closed translates into the following property of $\Delta$.
\par First observe that, because of transitivity of the action, we have $\Delta=\{ 0_\Delta+\gamma : \gamma\in \Gamma\}$. So, given $\delta\in \Delta$, let $\delta:=0_{\Delta}+\gamma_{\delta}$ for some $\gamma_{\delta}\in \Gamma$.
Since $\Gamma$ is o.l.-closed, there exists $\gamma_{\mu}\in \Gamma$ such that $\gamma_{\delta}=\min_{0\leq i\leq d}\sigma_{v}^i(\gamma_{\mu})+\gamma_{i}.$ Let $0\leq i_{0}\leq d$ be such that $\sigma_{v}^{i_{0}}(\gamma_{\mu})+\gamma_{i_{0}}=\min_{0\leq i\leq d}\sigma_{v}^i(\gamma_{\mu})+\gamma_{i}.$
Let $\mu:=0_{\Delta}+\gamma_{\mu}$; we have $\delta=0_{\Delta}+\gamma_{\delta}=0_{\Delta}+\sigma_{v}^{i_{0}}(\gamma_{\mu})+\gamma_{i_{0}}=\tau^{i_{0}}(\mu)+\gamma_{i_{0}}$ (cf. Definition \ref{langageLV}). 
  Since the action of $\Gamma$ on $\Delta$ respects the order, we have also that $\delta=\min_{0\leq i\leq d}\tau^i(\mu)+\gamma_{i}.$

\defn We will say that $(\Delta,\le, 0_{\Delta},\tau,+\gamma;\gamma\in \Gamma)$ is ordered linearly closed (o.l.-closed) if 
given any finite subset $\{\gamma_{i}\in\Gamma; 0\leq i\leq d\}$,  for any $\delta\in \Delta$ there exists $\mu\in \Delta$ such that $\delta=\min_{0\leq i\leq d}\{\tau^i(\mu)+\gamma_{i}\}.$ 
\edefn
\nota \label{notation-upsilon}
\par Given $q(t)\in A$, $q(t)=\sum_{i=0}^d t^i a_{i}$, given $\delta$ and $\{\gamma_{i}:=v(a_{i}): 0\leq i\leq d\}$, we will denote the corresponding $\mu$ with  $\delta=\min\{ \tau^i(\mu)+\gamma_{i}:\;0\leq i\leq d\}$ by $\U(q,\delta)$. We also set $\U^{-1}(q,\mu):=\delta$.
\enota
 As soon as $\Delta$ is o.l.-closed, the functions $\U^{-1}$ and $\U$ are well-defined and we have the following relationship between $\U^{-1}$ and $\U$. Let $q(t)\in \I$, $\mu,\delta\in \Delta$, then
$\U^{-1}(q(t),\U(q(t),\delta))=\delta$ and $\U(q(t),\U^{-1}(q(t),\mu))=\mu$.
 Moreover, for $m\in M$ we always have $w(m.q(t))\geq \U^{-1}(q(t),w(m))$. And finally since each of the functions $\tau^i+\gamma_{i}$ are increasing on $\Delta$, $\U(q,\mu_{2})\leq \mu_{1}\leftrightarrow \mu_{2}\leq \U^{-1}(q,\mu_{1})$. This last equivalence implies that $\U$ is increasing and since it is injective, it is strictly increasing.
\bigskip
\par Note that if $\sigma_{v}$ is surjective on $\Gamma$, then because of the transitivity of the action, $\tau$ is surjective on $\Delta$. In the following lemma, we will show that $\tau$ surjective implies that $\Delta$ is o.l.-closed.
\lmm\label{angularbis} 
Assume that $\tau$ is surjective on $\Delta$. Then $\Delta$ is o.l.-closed.
So, for any $q(t)\in A$, the function $\U(q,.)$ is well-defined on $\Delta$ 
(and so strictly increasing).

\elmm
 \pr  Let $q(t)\in A$, $q(t)=\sum_{i=0}^d t^i a_{i}$, and let $\{\gamma_{i}:=v(a_{i}): 0\leq i\leq d\}$.
 Given $\delta_{0}\in \Delta$, let us show that there exists $\delta\in \Delta$ such that $\delta_{0}=\min_{0\le i\le d} \, \tau^i(\delta)+\gamma_{i}$.
\par Consider the functions $f_{i}$ on $\Delta$ defined by $f_{i}(\delta)=\tau^i(\delta)+\gamma_{i}$, $0\leq i\leq d$. Since $\tau$ is assumed here to be surjective, so is $\tau^i$. Thus there exists $\delta_{i}\in \Delta$ such that $\delta_{0}+(-\gamma_{i})=\tau^i(\delta_{i})$ and so
$\delta_{0}=f_{i}(\delta_{i})$. 
 \par Each function $f_{i}$ is strictly increasing: if $\delta_{1}<\delta_{2}$, then $\tau^i(\delta_{1})<\tau^i(\delta_{2})$ and $\tau^i(\delta_{1})+\gamma_{i}<\tau^i(\delta_{2})+\gamma_{i}$. So the maximum $\mu$ of the $\delta_{i}$'s such that $f_{i}(\delta_{i})=\delta_{0}$ is well-defined. 
 Since $\delta_{i}\leq \mu$, we have that $\delta_{0}=f_{i}(\delta_{i})\leq f_{i}(\mu)$ ($f_{i}$ is increasing), $0\leq i\leq d$. So, $\delta_{0}=\min_{0\leq i\leq n} f_{i}(\mu)$, namely $\U(q,\delta_{0})=\mu$.



\qed

\bigskip
\bigskip
\par In the following, we will examine the case when for $n\in M$ with $w(n)=\delta$, there exists $m\in M$ with $n=m\cdot q(t)$ and $w(m)=\U(q(t),\delta)$.


\defn Assume that $\tau$ is surjective on $\Delta$ and that $(\Delta,\le, 0_{\Delta},\tau,+\gamma;\gamma\in \Gamma)\models T_{\Delta,\tau}$, 
then let 
$T_V^+$ be the following $\L_V$-theory:
\begin{enumerate}
\item $T_{V}$,
\item\label{D} $\forall n\;(n\in V_{0}\rightarrow (\exists m\;(m\in V_{0}\,\&\,m\cdot q(t)=n)))$, for all $q(t)\in \I$, $q(t)$ $\sigma$-separable.


\end{enumerate}
\edefn
\lmm\label{wdiv}  Let $M$ be a valued $A$-module with a $\sigma$-basis, and which is a model of $T_V^+$. 
Let $q(t)\in \I$ be $\sigma$-separable and $\delta\in \Delta$. Let $\mu:=\U(q,\delta)$. Then, for any $n\in M$ with $w(n)=\delta$, there exists $m\in M$ such that $m\cdot q(t)=n\;\&\;w(m)=\mu$. Moreover $\mu$ has the additional property that for any $m\in M$ with $w(m)=\mu$, $w(m\cdot q(t))\geq \delta$. 
\elmm
 \pr Let $q(t)=\sum_{i=0}^d t^i a_{i}\in \I$, with $a_{0}\neq 0$. By axiom scheme (2) of $T_V^+$, 
 for any $n_{0}\in M$ with $w(n_{0})=0$, there exists $m_{0}\in M$ such that $m_{0}\cdot q(t)=n_{0}$ and $w(m_{0})=0$. 
 Let $k_{\delta}, k_{\mu}\in K$ be such that $0+v(k_{\delta})=\delta$. and $0+v(k_{\mu})=\mu$.
 \par Now let $n\in M$ with $w(n)=\delta$ and consider the polynomial $\tilde q(t):=\sum_{i} t^i k_{\mu}^{\sigma^i} a_{i} k_{\delta}^{-1}$. Then by construction $\tilde q(t)\in \I$ and is still $\sigma$-separable. Now $w(n\cdot k_{\delta}^{-1})=0$, so by hypothesis there exists $m_{0}\in M_{0}$ such that $m_{0}\cdot \tilde q(t)= n\cdot k_{\delta}^{-1}$. So, $m_{0}\cdot \sum_{i} t^i k_{\mu}^{\sigma^i} a_{i}=n$ and so $m_{0}\cdot k_{\mu} q(t)=n$. Set $m:=m_{0}\cdot k_{\mu}$, we have $w(m)=\mu$ and $m\cdot q(t)=n$.
 \par Moreover, if $w(m)=\mu$, then by the compatibility condition between $\tau$ and $\sigma_{v}$, we have $w(m\cdot q(t))\geq \min_{i} \tau^i(w(m))+v(a_{i})\geq \delta$.
 \qed
\kor\label{div}  Let $M$ be a model of $T_V^+$ 
as in the lemma.
Let $q(t)\in \I$ be $\sigma$-separable and $\delta\in \Delta$. 
Then, for any $n\in M_{\delta}$, there exists $m\in M_{\U(q,\delta)}$ such that $m\cdot q(t)=n$. Moreover $\U(q,\delta)$ is such that for any $m\in M_{\U(q,\delta)}$, $m\cdot q(t)\in M_{\delta}$. \qed
\ekor
\medskip

\par The separable closure of a valued field of characteristic $p$ is dense in its algebraic closure. This translates as follows in the case of models of $T_V^+$.

\lmm\label{actiontw} If 
$M$ is a valued $A$-module  which is a model of $T_V^+$, and if the action of $\Gamma$ on $\Delta$ satisfies the following $(\star\star)$: for all $\delta_{1}, \delta_2\in \Delta$, there exists $\gamma\in \Gamma$ such that  
$\delta_{1}+(\sigma_{v}(\gamma)-\gamma)=\delta_{2}.$ 
Then for any $\delta$ and  $m$ with $w(m)\leq \delta$, there exists $n$ such that $w(m-n\cdot t)=\delta.$ 
\elmm
\pr W.l.o.g. we may assume that $\delta>w(m)$ (otherwise it suffices to choose $n=0$). First choose 
$k\in K$ such that $w(m)+v(k^{\sigma})<0$.
\par By Lemma \ref{wdiv} $M$ is $(t-1)$-divisible, there exists $n\in M$ such that $m\cdot k^{\sigma}=n\cdot (t-1)$ and necessarily $w(n)<0$. So $\tau(w(n))<w(n)$ and therefore $w(n\cdot t)=w(m)+\sigma_{v}(v(k))$ and $w(n)=\tau^{-1}(w(m)+\sigma_{v}(v(k_{\mu})))$.
\par We have that $m=n\cdot k^{-1} t-n\cdot k^{-\sigma}$.  
\par We have $w(m-n\cdot k^{-1} t)=w(n)-v(k^{\sigma})=\tau^{-1}(w(m))+(v(k)-\sigma_{v}(v(k)))$. 
Now, by the extra assumption, there is  $k_1$ such that $\tau^{-1}(w(m))+(v(k_1)-\sigma_{v}(v(k_1)))=\delta$ 
It suffices to see that this forces $w(m)+v(k_1^\sigma)<0$, so that the preceding discussion applies to $k_1$ as well, and we are done. But we have $w(m)+v(k_1^\sigma) =\tau (\delta+v(k_1^\sigma))$ and $w(m)+v(k_1^\sigma)<\delta+v(k_1^\sigma)$. So if $w(m)+v(k_1^\sigma) \ge 0$, we would get $0<\delta+v(k_1^\sigma)<\tau(\delta+v(k_1^\sigma))$, and then  $w(m)+v(k_1^\sigma) <\tau (\delta+v(k_1^\sigma))$ which is absurd. Hence we must have $w(m)+v(k_1^\sigma)<0$, as wanted.
\qed
\bem Note that the action of $\Gamma$ on $\Delta$ is transitive, so to meet the hypothesis $(\star\star)$,
we can require that the action of $\sigma_{v}-1$ is surjective on $\Gamma$. 
\ebem
\defn Recall that $\Gamma^+:=\{\gamma\in \Gamma: \gamma\geq 0\}$. We will say that $\si_{v}$ is $2$-contracting on $\Gamma$ if $\forall \gamma\in \Gamma^+\;\si_{v}(\gamma)\geq \gamma+\gamma$.
\edefn
\kor\label{actiont} If $M\models T_V^+$ 
 and if $\si_{v}$ is $2$-contracting on $\Gamma$,
then for any $\delta\in \Delta$ and for all $m\in M$, there exists $n$ such that $\;V_{\delta}(m-n\cdot t)$ holds. 
\ekor
\pr It follows from the proof of the above lemma, noting that we only need in this setting that $\forall \delta_{1}\in \Delta\forall \delta_{2}\in \Delta\exists\;\gamma\in \Gamma\;\delta_{1}+(\sigma_{v}(\gamma)-\gamma)\geq\delta_{2}.$ W.l.o.g., we may assume that $\delta_{1}<\delta_{2}$. So given $\delta_{1}<\delta_{2}\in \Delta$, by transitivity of the action of $\Gamma$ on $\Delta$, we get that there exists $\tilde \gamma\in \Gamma^+$ such that $\delta_{1}+\tilde \gamma=\delta_{2}$. Since $\si_{v}$ is $2$-contracting we get that $\si_{v}(\tilde \gamma)\geq \tilde \gamma+\tilde \gamma$. So $\delta_{1}+\si_{v}(\tilde \gamma)\geq\delta_{2}+\tilde \gamma$ and so $\delta_{1}+\si_{v}(\tilde \gamma)-\tilde \gamma\geq\delta_{2}.$


\bigskip
\par We now check that the basic example of the separable closure of $\tilde \F_{p}((T))$ is not covered by Rohwer (see \cite{R}, pp. 40-41), since it does not have a weak valuation basis. 
\defn\label{Tw,tau,alpha} 
Let $\M:=(M,+,0,.r; r\in A, \lambda_{i}, i\in n)$ be a model of $T_{\sigma}$.
Then $\M$ is a {\em valued $A$-module with a weak $\si$-valuation basis} 
if there exists $r\in K$ such that for each $m\in M$ we have: $w(m)\leq \min_{i} \{w(\lambda_{i}(m)\cdot t)+v(c_{i})+v(r):\;c_{i}\in \C\}.\;\;(*)$
\edefn
\lmm Let $K$ be any valued separably closed field $K$ of finite imperfection degree, then $K$ does not have a weak $\si$-valuation basis, with $\si$ the Frobenius endomorphism.
\elmm
\pr By way of contradiction, let $c_1, c_2, \ldots$  be a linear basis of $K$ over $K^p$ and suppose that it is a weak $\si$-valuation basis
and let $\delta$ be the corresponding $v(r)$. By adjusting $\delta$ and since $v(K^{\times})$ is $p$-divisible, we may always assume that $c_{1}=1$. 
Let $\delta'\in \Gamma$ such that $\delta' > \{v(c_2), \delta, v(c_2)+\delta\}$.  By Corollary \ref{actiont} let $a, b$ such that $c_2=a^p+b$ with $v(b)\ge \delta'$. If $v(a^p)\ne v(c_2)$, then that would contradict the required inequality (*) for $a^p-c_2$. Otherwise, $v(a^p)=v(c_2)$, and again this contradicts (*) for $a^p-c_2$.\qed


\defn Assume that $\tau$ is surjective on $\Delta$, and that $\si_{v}$ is $2$-contracting on $\Gamma$. Then we let $T_{V}^{sep}$ be the $\L_V$-theory $T_V^+$.
\edefn

In particular, if $K$ is separably closed of characteristic $p$ and $\sigma$ is the  Frobenius map, then $M=K$ is a model of $T_V^{sep}.$



We will prove in the next sections that $T_V^{sep}$ eliminates quantifiers up to index sentences. 




\section{Special cases}

\par In order to eliminate quantifiers in $T_V^{sep}$, we need some basic cases and reductions, which are treated in the following lemmas. 
Our main tools will be Corollaries \ref{div} and \ref{actiont} and we will use Notation \ref{notation-upsilon}. We will treat the general case in the next section. 

\par We will use the notation
$u\cdot r\equiv_\delta m$ to mean
that 
$V_{\delta}(u\cdot r-m)$ holds.

\lmm\label{cas1} Consider a system of the form
\begin{equation*}
\exists u\left\{
\begin{array}{rl}
u\cdot t\equiv_{\mu_{1}} b_{1}\\
u\cdot r\equiv_{\mu_{2}} b_{2},\\
\end{array}\right.
\end{equation*}
where $r\in \I$ is separable. Then this system is equivalent to a congruence of the form $$b_1\cdot r^\sigma \equiv_{\mu_3} b_2\cdot t.$$
\elmm

\pr We distinguish two cases.
\par $(i)$ $\U^{-1}(r^{\sigma},\mu_{1})\geq \tau(\mu_{2})$.  Then the system above is equivalent to $$b_{1}\cdot r^{\sigma}\equiv_{\tau(\mu_{2})}b_{2}\cdot t.$$
 One implication is straightforward. For the reverse implication, 
since $M\models T_{V}^{sep}$, by Corollary \ref{actiont}, there exists $u$ such that $\,u\cdot t\equiv_{\mu_{1}}b_{1}$. So $u\cdot t r^{\sigma}\equiv_{\U^{-1}(r^{\sigma},\mu_{1})}b_{1}\cdot r^{\sigma}$. So $u\cdot r t\equiv_{\tau(\mu_{2})}b_{2}\cdot t$ and so $u\cdot r\equiv_{\mu_{2}}b_{2}.$

\par $(ii)$ $\U^{-1}(r^{\sigma},\mu_{1})< \tau(\mu_{2}).$ Then the system above is equivalent to $$b_{1}\cdot r^{\sigma}\equiv_{\U^{-1}(r^{\sigma},\mu_{1})}b_{2}\cdot t.$$
 One implication is straighforward. For the reverse implication, choose $\mu$ such that $\U^{-1}(r^{\sigma},\mu)\geq \tau(\mu_{2})$ (and so $\mu>\mu_{1}$). Again by Corollary \ref{actiont}, there exists $u$ such that $u\cdot t\equiv_{\mu_{1}}b_{1}$. So $u\cdot t r^{\sigma}\equiv_{\U^{-1}(r^{\sigma},\mu_{1})} b_{1}\cdot r^{\sigma}\equiv_{\U^{-1}(r^{\sigma},\mu_{1})} b_{2}\cdot t$.
\par Since $M\models T^{sep}$, by Corollary \ref{div}, there exists $u''$ such that $u''\cdot r^{\sigma}=(u\cdot r-b_{2})\cdot t$ with $w(u'')\geq \mu_{1}$ (by definition of $\U^{-1}(r^{\sigma},\mu_{1})$). 
 Let $u'$ be such that $u'\cdot t\equiv_{\mu} u''$.
Then $u'\cdot r t=u'\cdot t r^{\sigma}\equiv_{\tau(\mu_{2})} (u\cdot r-b_{2})\cdot t$, which implies that $(u-u')\cdot r t\equiv_{\tau(\mu_{2})} b_{2}\cdot t$ and so $(u-u')\cdot r\equiv_{\mu_{2}} b_{2},$ which finishes the proof since $(u-u')\cdot t\equiv_{\mu} u\cdot t-u''\equiv_{\mu_{1}} b_{1}$.
\qed

\lmm\label{cas2} Consider a system of the form
\begin{equation*}
\exists u\left\{
\begin{array}{rl}
u\cdot r_{1}= b_{1}\\
u\cdot r_{2}\equiv_{\delta_{2}} b_{2},\\
\end{array}\right.
\end{equation*}
where $r_{1},\;r_{2}\in \I$ are separable and assume that $deg(r_{1})\geq deg(r_{2})$.
 Then this system is equivalent to the following system 
 \begin{equation*}
\exists u\left\{
\begin{array}{rl}
u\cdot r_{2}= b_{2}\\
u\cdot r_{3}\equiv_{\delta} b_{1}\cdot \lambda-b_{2}\cdot s,\\
\end{array}\right.
\end{equation*}
where $\delta=\U^{-1}(r_{1}\lambda,\U(r_{2},\delta_{2}))$, for some $\lambda\in {\mathcal O}_K$, and $deg(r_2)> deg(r_3)$.
\elmm

\pr  By the generalized euclidean algorithm, there exists $\lambda\in \O_{K}$ such that $r_{1}\lambda=r_{2} s+r_{3}$ with $deg(r_{3})<deg(r_{2})$.  
\par Suppose $u$ is a solution of the first system. Let $u'$ be such that $u'\cdot r_{2}=u\cdot r_{2}-b_{2}$. We can find such $u'$ with $w(u')\geq \U(r_{2},\delta_{2})$.
So $(u-u')\cdot r_{2}.s+(u-u')\cdot r_{3}=(u-u')\cdot r_{1} \lambda=b_{1}\cdot \lambda-u'\cdot r_{1} \lambda$ i.e. $(u-u')\cdot r_{3}=b_{1}\cdot \lambda-b_{2}\cdot s-u'\cdot r_{1} \lambda$.
\par Conversely, let $u''$ satisfy the second system. Then $u''\cdot r_{1} \lambda=u''\cdot r_{2} s+u''\cdot r_{3}\equiv_{\delta} b_{1}\cdot \lambda$.
So let $u'''$ be such that $u'''\cdot r_{1}=u''\cdot r_{1}-b_{1}$ and we have to make sure that we can choose $u'''$ such that $w(u'''\cdot r_{2})\geq \delta_{2}$. In other words, $\U(r_{1} \lambda,\delta)= \U(r_{2},\delta_{2})$. \qed


\lmm\label{cas3} Consider a system of the form
\begin{equation*}
\exists u\left\{
\begin{array}{rl}
u\cdot r_{1}\equiv_{\delta_{1}} b_{1}\\
u\cdot r_{2}\equiv_{\delta_{2}} b_{2}\\
\end{array}\right.
\end{equation*}
where $r_{1},\;r_{2}\in \I$ are separable and $\U(r_{1},\delta_{1})\leq \U(r_{2},\delta_{2})$.
 Then this system is equivalent to the following system 
\begin{equation*}
\exists u\left\{
\begin{array}{rl}
u\cdot r_{1}\equiv_{\delta_{1}} b_{1}\\
u\cdot r_{2}= b_{2}\\
\end{array}\right.
\end{equation*}
\elmm

\pr Indeed, we can choose $u'$ such that $w(u')\geq \U(r_{2},\delta_{2})$ and such that $u'\cdot r_{2}=b_{2}-u\cdot r_{2}$ and so $(u+u')\cdot r_{2}=b_{2}$. Moreover $w(u'\cdot r_{1})\geq \delta_{1}$. \qed

\lmm\label{cas6} Consider a system of the form
\begin{equation*}
\exists u\left\{
\begin{array}{rl}
u\cdot r= b\\
u\cdot t\equiv_{\delta} b_{1},\\
\end{array}\right.
\end{equation*}
 where $r\in \I$ is separable. 
 Then this system is equivalent to the following system 
\begin{equation*}
\exists u\left\{
\begin{array}{rl}
u\cdot r\equiv_{\delta'} b\\
u\cdot t\equiv_{\delta} b_{1},\\
\end{array}\right.
\end{equation*}
where $\delta'$ is chosen such that $\tau(\U(r,\delta'))\geq \delta$.
\elmm

\pr Let us show the non-trivial implication. Since $r$ is separable, there exists $u'$ such that $u'\cdot r=(u\cdot r-b)$ and we can choose such $u'$ with $w(u')\geq \U(r,\delta')$. Since $\tau(\U(r,\delta'))\geq \delta$, we get that $w(u'\cdot t)\geq \delta$ and so $(u-u')$ is a solution of the first system. \qed



\lmm\label{cas_constante} Consider a system of the form
 \begin{equation*}
\label{6}\exists u\left\{
\begin{array}{rl}
u\cdot r_{0}= b_{0}\\
\bigwedge_{i=1}^d u\cdot t^{n}\equiv_{\delta_{i}} b_{i},\\
\end{array}\right.
\end{equation*}
where all $r_{i}\in \O_{K}$, $\bar r_i \ne 0$, $1\leq i\leq d$. Then this system is equivalent to congruences of the following form
\begin{equation*}
\left\{
\begin{array}{rl}
b'_{0}\cdot t^{n}\equiv_\delta b'_{1}\cdot r_{0}^{\sigma^n}\\
\bigwedge_{i=2}^d b'_{1}\equiv_{\delta'_{i}} b'_{i},\\
\end{array}\right.
\end{equation*}
where the $b'_i$ have the same parameters as the $b_i$.
\elmm

\pr We first proceed as in Lemma \ref{cas6}, replacing the equation 
$u\cdot r_{0}=b_{0}$ by a congruence $u\cdot r_{0}\equiv_{\delta'} b_{0}$ where $\delta'$ is chosen such that $\tau^n(\Upsilon(r_{0},\delta'))\geq\{\delta_{i}: 1\leq i\leq d\}$.
\par So it remains to consider a system of the form:
\begin{equation*}
\label{7}\exists u\left\{
\begin{array}{rl}
u\cdot r_{0}\equiv_{\delta'} b_{0}\\
\bigwedge_{i=1}^d u\cdot t^{n}\equiv_{\delta_{i}} b_{i},\\
\end{array}\right.
\end{equation*}
Now we proceed as in Lemma \ref{cas1}. First we note that $u\cdot r_{0} t^n\equiv_{\tau^n(\delta')} b_{0}\cdot t^n$ is equivalent to $u\cdot r_{0}\equiv_{\delta'} b_{0}$. We rewrite the first formula as $u\cdot t^n r_{0}^{\sigma^n}\equiv_{\tau^n(\delta')} b_{0}\cdot t^n$.
\par We order the $\delta_{i}$ and w.l.o.g. assume that $\delta_{1}\geq\max\{\delta_{i}: 1\leq i\leq d\}$.
Our system 
is then equivalent to:
\begin{equation*}
\label{8}\left\{
\begin{array}{rl}
b_{0}\cdot t^{n}\equiv_{\Upsilon^{-1}(r_{0}^{\sigma^n},\delta_{1})} b_{1}\cdot r_{0}^{\sigma^n}\\
\bigwedge_{i=2}^d b_{1}\equiv_{\delta_{i}} b_{i},\\
\end{array}\right.
\end{equation*}
Indeed, by Corollary \ref{actiont}, there exists $u$ such that $u\cdot t^n\equiv_{\delta_{1}} b_{1}$. So, if $\Upsilon^{-1}(r_{0}^{\sigma^n},\delta_{1})\geq \tau^n(\delta')$ $(\star)$, we get that $u\cdot t^n r_{0}^{\sigma^n}\equiv_{\Upsilon^{-1}(r_{0}^{\sigma^n},\delta_{1})} b_{1}\cdot r_{0}^{\sigma^n}$ and so 
since $b_{1}\cdot r_{0}^{\sigma^n}\equiv_{\Upsilon^{-1}(r_{0}^{\sigma^n},\delta_{1})}b_{0}\cdot t^n$, we get by $(\star)$, that $u\cdot r_{0}\equiv_{\delta'} b_{0}.$
\par Now assume that $\Upsilon^{-1}(r_{0}^{\sigma^n},\delta_{1})< \tau^n(\delta')$. So we choose $\delta''$ such that $\Upsilon^{-1}(r_{0}^{\sigma^n},\delta'')\geq \tau^n(\delta')$.
\par Again, by Corollary \ref{actiont}, there exists $u$ such that
$u\cdot t^{n}\equiv_{\delta_{1}}b_{1}$. So we get $$u\cdot t^n r^{\sigma^n}\equiv_{\U^{-1}(r_{0}^{\sigma^n},\delta_{1})}b_{1}\cdot r^{\sigma^n}\equiv_{\U^{-1}(r_{0}^{\sigma^n},\delta_{1})} b_{0}\cdot t^n.$$
By Corollary \ref{div}, there exists $u''$ such that $u''\cdot r^{\sigma^n}=(u\cdot r-b_{0})\cdot t^n$ with $w(u'')\geq \delta_{1}$. 
\par By Corollary \ref{actiont}, there exists $u'$ such that $u'\cdot t^n\equiv_{\delta''} u''$.
Then $u'\cdot r t^n=u'\cdot t^n r^{\sigma^n}\equiv_{\tau^n(\delta')} (u\cdot r-b_{0})\cdot t^n$, which implies that $(u-u')\cdot r t^n\equiv_{\tau^n(\delta')} b_{0}\cdot t^n$ and so $(u-u')\cdot r\equiv_{\delta'} b_{0}.$
\par Since $w(u'')=\delta_{1}$, we may add to the other congruences $u''=u'\cdot t^n$ without perturbing them. \qed

\section{Quantifier elimination}

\par We now prove that $T_{V}^{sep}$ admits quantifier-elimination up to index sentences. 

\nota\label{lambda2} Let $\boldsymbol d=(d_{1},\cdots,d_{m})\in n^m$ be a $m$-tuple of natural numbers between $0$ and $n-1$. We denote by $\lambda_{\boldsymbol d}^{(m)}$ the composition of the $m$ $\lambda$-functions: $\lambda_{d_{1}}\circ\lambda_{d_{2}}\circ\cdots\circ\lambda_{d_{m}}$.
\enota
\lmm\label{sep} A system of equations $\bigwedge_{i=1}^d u\cdot r_{i}=t_{i}(\bar y)$, where $t_{i}(\bar y)$ is a $\L_{A}$-term and $r_{i}\in A$ with at least one $r_{i}$ $\sigma$-separable, is equivalent to one equation of the form $u\cdot r=t(\bar y)$, where $r\in A$ is separable together with a conjunction of atomic formulas in $\bar y$. 
\elmm
\pr We apply the Euclidean algorithm and do some bookeeping to check that we always keep a separable coefficient. Assume that $r_{1}$ is separable.
Let us consider the system:
\begin{equation}
\left\{
\begin{array}{rl}
u\cdot r_{1}= t_{1}(\bar y)\\
u\cdot r_{i}= t_{i}(\bar y),\\
\end{array}\right.
\end{equation}
with $i\neq 1$.
\par If $r_{i}$ is not separable and if $deg(r_{1})\geq deg(r_{i})$, then for some $r',\;r''\in A$, we have $r_{1}=r_{i} r'+r''$, then $r''\neq 0$ and $r''$ is separable. 
So, the system is equivalent to:
\begin{equation}
\left\{
\begin{array}{rl}
u\cdot r_{i}= t_{i}(\bar y)\\
u\cdot r''=t_{1}(\bar y)-t_{i}(\bar y)\cdot r',\\
\end{array}\right.
\end{equation}
with $deg(r'')<deg(r_{i})$ and $r''$ separable.
\par If $r_{i}$ is not separable and if $deg(r_{1})< deg(r_{i})$, then for some $r',\;r''\in A$, we have $r_{i}=r_{1} s'+s''$, then either $s''= 0$ and
the system is equivalent to:
\begin{equation}
\left\{
\begin{array}{rl}
u\cdot r_{1}= t_{1}(\bar y)\\
t_{1}(\bar y)\cdot s'=t_{i}(\bar y),\\
\end{array}\right.
\end{equation}
or $s''\neq 0$ and the system is equivalent to:
\begin{equation}
\left\{
\begin{array}{rl}
u\cdot r_{1}= t_{1}(\bar y)\\
u\cdot s''=t_{i}(\bar y)-t_{1}(\bar y)\cdot s'.\\
\end{array}\right.
\end{equation}
\par If $r_{i}$ is separable, then w.l.o.g. $deg(r_{1})\geq deg(r_{i})$. For some $r',\;r''\in A$, we have $r_{1}=r_{i} r'+r''$. Either $r''=0$ and the system is equivalent to:
\begin{equation}
\left\{
\begin{array}{rl}
u\cdot r_{i}= t_{i}(\bar y)\\
t_{1}(\bar y)=t_{i}(\bar y)\cdot r',\\
\end{array}\right.
\end{equation}
or $r''\neq 0$ and the system is equivalent to 
\begin{equation}
\left\{
\begin{array}{rl}
u\cdot r_{i}= t_{i}(\bar y)\\
u\cdot r''=t_{1}(\bar y)-t_{i}(\bar y)\cdot r',\\
\end{array}\right.
\end{equation}
with $r_{i}$ separable.
\par In each case, we showed that the system of two equations with the pair of coefficients $(r_{1},s_{i})$ where $r_{1}$ separable, was equivalent with another system with a pair of coefficients consisting of a separable coefficient and such that the sum of the degrees of the coefficients decreased. If one of the coefficient is zero, we consider another equation, if applicable, of the conjunction. If both coefficients are nonzero, we repeat the procedure until either we considered all of the equations occurring in the conjunction, or one of the coefficient has degree zero which allows us to eliminate the variable $u$.
\qed
\prop \label{EQsep} In $T_{V}^{sep}$, every $\L_{V}$ p.p. formula is equivalent to a positive quantifier-free formula.
\eprop
\pr As usual, we proceed by induction on the
number of existential quantifiers, so it suffices
to consider a formula existential in just one variable
$\exists u \phi(u,\bold y)$, where $\phi(u,\bold y)$ is a
conjunction of atomic $\L_{V}$-formulas. 
\par Note first that terms in $u$ are $L_{A}$-terms in $u$, $\lambda_{i}(u)$, $i\in n^{\ell}$, for some $\ell\geq 1$, where $\lambda_{i}$ denotes the composition of $\ell$ functions $\lambda_{j}$, $j\in n$ (see \cite[Notation 3.3]{DDP1}). One uses the fact that the $\lambda_{i}$ functions are additive and that $\lambda_{i}(u\cdot q(t))$ with $q(t)\in A$, can be expressed as an $L_{A}$-term in $\lambda_{j}(u)$, $j\in n$. Moreover since $u=\sum_{i\in n}\lambda_{i}(u)\cdot t c_{i}$, we may assume that the terms are terms in only the $\lambda_{i}(u)$, $i\in n^{\ell}$ (see \cite{DDP1} Lemma 3.2,  and Notation \ref{lambda}). 
\par Therefore we may replace the quantifier $\exists u$ by $n^{\ell}$ quantifiers $\exists u_{n^{\ell-1}}\cdots\exists u_{0}\, \bigwedge_{i\in n^{\ell}} u_{i}=\lambda_{i}(u)$. We first tackle the quantifier $\exists u_{0}$ and for convenience, let us replace $u_{0}$ by $u$.
\par Since $A$ is right Euclidean, we can always
assume that we have at most one atomic formula involving
$u_{0},$ of the form $u_{0}\cdot r_{0}(t)=t_{0}(\bold{y})$, where $t_{0}(\bold{y})$ is a $\L_{A}$-term. 
\cl\label{sepequation}  We may assume that $r_{0}$ is separable. 
\ecl
\prcl Write $r_{0}=t^m r_{0}'$, where $m\in \N$ and $r_{0}'$ separable.
Express $r_{0}'=\sum_{\boldsymbol d\in n^m} r_{0\boldsymbol d}' c_{\boldsymbol d}$ with the property
that $r_{0\boldsymbol d}' \in A^{\sigma^m}[t;\sigma]$ e.g. $r_{0\boldsymbol d}'=\sum_j t^j a_{\boldsymbol d j}^{\sigma^m}$,
with $a_{\boldsymbol d j} \in K$. 
Recall that $\sqrt[\sigma^m]{r_{0\boldsymbol d}'}=\sum_j t^j a_{\boldsymbol d j}$, so
$t^m r_{0}'=\sum_{\boldsymbol d\in n^m} \sqrt[\sigma^m]{r_{0\boldsymbol d}'} t^m c_{\boldsymbol d}$ (see Notation \ref{lambda}).
\par Using this equality, replace the atomic formula $u\cdot t^m r_{0}'=t_{0}$ by the system $$\bigwedge_{\boldsymbol d\in n^m}\;u\cdot \sqrt[\sigma^m]{r_{0\boldsymbol d}'}=\lambda_{\boldsymbol d}^{(m)}(t_{0})$$ (see Notation \ref{lambda2}). Note that for at least one tuple $\boldsymbol d$, $\sqrt[\sigma^m]{r_{0\boldsymbol d}'}$ is separable. So by Lemma \ref{sep}, we may assume that we have just one equation with a separable coefficient together a conjunction of atomic formulas in $\bar y$. \qed
\medskip
\par Moreover, for any element $r(t)=\sum_{j} t^j a_{j}\in A$, there exists $\mu\in  K$ such that $r(t).\mu\in \I$ (let $\mu:=a_{k}^{-1}$ where $v(a_{k})=\tilde v(r(t))$ and $\I$ denote the set of elements of $A_{0}$ which are non trivial residually, in other words which have a coefficient with valuation zero. So, we transform the atomic formula: $u\cdot r(t)=t(\bold{y})$ by multiplying both sides by $\mu$ and $V_{\delta}(u\cdot r-
t(\bold{y}))$ by $V_{\delta+v(\mu)}(u\cdot r(t) \mu-
t(\bold{y})\cdot \mu),$ $\delta\in \Delta$.
\par So we reduced ourselves to consider an existential formula of the form 
$\exists u\;\phi(u,\bold{y})$, where $\phi(u,\bold{y})$ is of the form
$$u\cdot r_{0}=t_{0}(\bold{y})\;\&\; \bigwedge_{k=1}^m\;V_{\delta_{k}}(u\cdot r_{k}-
t_{k}(\bold{y}))\;\&\;
\theta(\bold{y})$$
with $r_{k}\in \I$, $\theta(\bold{y})$ a
quantifier-free $\L_{V}$-formula,
$t_{k}(\bold{y})$ are $\L_{A}$-terms, and
$\delta_{k}\in \Delta$.


\par Note also that in case $r_{0}\neq 0$, we can always assume that
$deg(r_0)> deg(r_k)$, for all $k$. Indeed, suppose 
that
$deg(r_{0})\leq deg(r_{k}),$ for some $k$, say $k=1$.
By the g.r.
Euclidean algorithm in $A_0$, there exists
$\mu\in\O_{ K}$ such that
$r_{1} \mu=r_{0} r+r_{1}'$ with
$deg(r_{1}')<deg(r_{0})$ and $r,\;r_{1}'\in A_{0}$.
So, we have that
$u\cdot r_{1} \mu=u\cdot r_{0} r+u\cdot r_{1}'=t_{0}\cdot r+u\cdot r_{1}'$,
and we can replace $V_{\delta_1}(u\cdot r_1-t_1)$ by  
$V_{\delta_{1}+v(\mu)}(u\cdot r_{1}'+t_{0}(\bold{y})\cdot r-
t_{1}(\bold{y})\cdot \mu)$.  
\par First, we will assume that the
equation present in $\phi(u,\bold y)$, $u.r_{0}=t_{0}(\bold{y})$, is non trivial, namely that $r_{0}\neq 0$. We will
concentrate on the system formed by this equation and 
the  {\it congruences}. For ease of notation, we replaced $t_{0}(\bold{y})$ by $b_{0}$ and $t_{i}(\bold{y})$ by $b_{i}$. So consider a system of the form 
\begin{equation}
\label{0}\exists u\left\{
\begin{array}{rl}
u\cdot r_{0}= b_{0}\\
\bigwedge_{i=1}^d u\cdot t^{n_{i}} r_{i}\equiv_{\delta_{i}} b_{i},\\
\end{array}\right.
\end{equation}
 where $r_{0},\;r_{i}\in \I$,  $r_{0},\;r_{i}$ are $\sigma$-separable, $n_{i}\in \N$, $1\leq i\leq d$. 
\par We will call $\sum_{i=0}^d deg(r_{i})$ the {\it separability degree} of that system, and we proceed by induction on that number.
\par We consider two cases : either there is $1\leq i\leq d$ such that $n_{i}\geq 1$, or for all $1\leq i\leq d$, $n_{i}=0$.  We will refer to the latter systems as {\it separable} systems, namely those for which $r_{0}, r_i$ are $\sigma$-separable and $n_{i}=0$ for all $1\leq i\leq n$.
\par {\bf Case A:}  let $n_{0}:=\max \{n_{i}: 1\leq i\leq d\}$ and suppose $n_{0}\geq 1$. Then there exists $\delta$ such that the system \eqref{0} is equivalent to 
\begin{equation}
\label{1}\exists u\left\{
\begin{array}{rl}
u\cdot r_{0}\equiv_{\delta} b_{0}\\
\bigwedge_{i=1}^d u\cdot t^{n_{i}} r_{i}\equiv_{\delta_{i}} b_{i},\\
\end{array}\right.
\end{equation}
where $\delta$ is chosen such that $\Upsilon(r_{0},\delta)\geq\max\{\tau^{-n_{i}}(\delta_{i}): 1\leq i\leq d\}$ (see Lemma \ref{cas6}). 
Then system \eqref{1} is equivalent to
\begin{equation}
\label{2}\exists u\left\{
\begin{array}{rl}
u\cdot r_{0} t^{n_{0}}\equiv_{\tau^{n_{0}}(\delta)} b_{0}\cdot t^{n_{0}}\\
\bigwedge_{i=1}^d u\cdot t^{n_{i}} r_{i} t^{n_{0}-n_{i}}\equiv_{\tau^{n_{0}-n_{i}}(\delta_{i})} b_{i}\cdot t^{n_{0}-n_{i}},\\
\end{array}\right.
\end{equation}
We re-write system \eqref{2} as follows:
\begin{equation}
\label{3}\exists u\left\{
\begin{array}{rl}
u\cdot t^{n_{0}} r_{0}^{\sigma^{n_{0}}}\equiv_{\tau^{n_{0}}(\delta)} b_{0}\cdot t^{n_{0}}\\
\bigwedge_{i=1}^d u\cdot t^{n_{0}} r_{i}^{\sigma^{n_{0}-n_{i}}}\equiv_{\tau^{n_{0}-n_{i}}(\delta_{i})} b_{i}\cdot t^{n_{0}-n_{i}},\\
\end{array}\right.
\end{equation}
 If all $r_{i}\in \O_{K}$, then we are done by Lemma \ref{cas_constante}.  Otherwise we replace $u\cdot t^{n_{0}}$ by $u_{0}$ and we consider the following separable system of congruences, assuming that one of the $r_{i}\notin \O_{K}$:
\begin{equation}
\label{4}\exists u_0\left\{
\begin{array}{rl}
u_{0}\cdot r_{0}^{\sigma^{n_{0}}}\equiv_{\tau^{n_{0}}(\delta)} b_{0}\cdot t^{n_{0}}\\
\bigwedge_{i=1}^d u_{0}\cdot r_{i}^{\sigma^{n_{0}-n_{i}}}\equiv_{\tau^{n_{0}-n_{i}}(\delta_{i})} b_{i}\cdot t^{n_{0}-n_{i}},\\
\end{array}\right.
\end{equation}
Suppose we can solve that system (see Lemma \ref{cas3})
. Then by Corollary \ref{actiont}, there exists $u$ such that $u\cdot t^{n_{0}}\equiv_{\delta_{0}}u_{0}$, where we can choose $\delta_{0}\geq \max_{1\leq i\leq d} \{\tau^{n_{0}}(\delta),\tau^{n_{0}-n_{i}}(\delta_{i})\}$.
\par Now we order the set $\{\Upsilon(r_{i}^{\sigma^{n_{0}-n_{i}}},\tau^{n_{0}-n_{i}}(\delta_{i})), \Upsilon(r_{0}^{\sigma^{n_{0}}},\tau^{n_{0}}(\delta)): 1\leq i\leq d\}$ and we replace
one of the congruences (the one corresponding to the maximum index) by the corresponding equation. Then we apply the g.r. Euclidean algorithm in order to obtain 
a system with separability degree strictly smaller than that of system (\ref{4}). 
\par {\bf Case B:} suppose that for all $i, n_i=0.$ We have a system of the form 
\begin{equation}
\label{5}\exists u\left\{
\begin{array}{rl}
u\cdot r_{0}= b_{0}\\
\bigwedge_{i=1}^d u\cdot r_{i}\equiv_{\delta_{i}} b_{i},\\
\end{array}\right.
\end{equation}
 We order the set $\{\Upsilon(r_{i},\delta_{i}): 1\leq i\leq d\}$ and we show that system \eqref{5} is equivalent to another system where both the degrees of $r_{0}$ and of $r_{i}$, $1\leq i\leq d$ have decreased.
 We proceed as in Lemma \ref{cas2} 
 with $\delta_{2}$ replaced by the element of $\Delta$ realizing the maximum of the set above. 
 Considering the equation with each of the congruences, we obtain a system where the coefficient appearing in the equation is separable and has degree strictly less than $r_{0}$ and the coefficients occurring in the congruences have each a strictly smaller degree but may no longer be separable. Then we use the previous Case A, to obtain an equivalent system where now all the coefficients of the congruences are separable and the degrees of the coefficients of both the equation and the congruences either stayed the same or have decreased. We obtain a system with strictly smaller separability degree.
\medskip
\par Second, we will assume that there is no
equation present in $\phi(u,\bold y)$. So, we consider
a system formed by 
{\it congruences} and for ease of notation, as before, we replace $t_{i}(\bold{y})$ by $b_{i}$. Consider a system of the form  
\begin{equation}
\label{9}\exists u\left.
\begin{array}{rl}
\bigwedge_{i=1}^d u\cdot t^{n_{i}} r_{i}\equiv_{\delta_{i}} b_{i},\\
\end{array}\right.
\end{equation}
 where $r_{i}\in \I$, $r_i$ is separable, $n_{i}\in \N$, $1\leq i\leq d$. 
 \par Again we distinguish the two cases : either there is $1\leq i\leq d$ such that $n_{i}\geq 1$, or for all $1\leq i\leq d$, $n_{i}=0$.
\par {\bf Case A':}  let $n_{0}:=\max \{n_{i}: 1\leq i\leq d\}$ and suppose $n_{0}\geq 1$. Then 
the system \eqref{9} is equivalent to 
\begin{equation}
\label{10}\exists u\left.
\begin{array}{rl}
\bigwedge_{i=1}^d u\cdot t^{n_{0}} r_{i} t^{n_{0}-n_{i}}\equiv_{\tau^{n_{0}-n_{i}}(\delta_{i})} b_{i}\cdot t^{n_{0}-n_{i}}\\
\end{array}\right.
\end{equation}
 First note that if all $r_{i}\in \O_{K}$, it implies since $r_{i}\in \I$, that $r_{i}^{-1}\in \O_{K}\cap \I.$ 
In this case, w.l.o.g. we may assume that system (\ref{9}) is of the form
 \begin{equation}
\label{11}\exists u\left.
\begin{array}{rl}
\bigwedge_{i=1}^d u\cdot t^{n_{i}}\equiv_{\delta_{i}} b_{i}\\
\end{array}\right.
\end{equation}
System (\ref{11}) is equivalent to the following system:
\begin{equation}
\label{12}\exists u\left.
\begin{array}{rl}
\bigwedge_{i=1}^d u\cdot t^{n_{0}}\equiv_{\tau^{n_{0}-n_{i}}(\delta_{i})} b_{i}\cdot t^{n_{0}-n_{i}}\\
\end{array}\right.
\end{equation}
We order the elements $\{\tau^{n_{0}-n_{i}}(\delta_{i}):\;1\leq i\leq n\}$. Let $\alpha$ be a permutation of $\{1,\cdots,n\}$ and suppose that $\tau^{n_{0}-n_{\alpha(1)}}(\delta_{\alpha(1)})\leq\cdots\leq \tau^{n_{0}-n_{\alpha(n)}}(\delta_{\alpha(n)})$.
We claim that system (\ref{12}) is equivalent to:
\begin{equation}
\label{13} 
\left.
\begin{array}{rl}
\bigwedge_{i=1}^{d-1} b_{\alpha(i)}\cdot t^{n_{0}-n_{\alpha(i)}}\equiv_{\tau^{n_{0}-n_{\alpha(i)}}(\delta_{\alpha(i)})} b_{\alpha(i+1)}.t^{n_{0}-n_{\alpha(i+1)}}\\
\end{array}\right.
\end{equation}
 We use Lemma \ref{actiont} in order to find $u$ such that $u\cdot t^{n_{0}}\equiv_{\tau^{n_{0}-n_{j}}(\delta_{j})} b_{j}\cdot t^{n_{0}-n_{j}}$ and then we use the congruences. 
\par Now assume that  $r_{i}\notin \O_{K}$ for some $i$ in system (\ref{10}). Replace $u\cdot t^{n_{0}}$ by $u_{0}$ and consider the system :
\begin{equation}
\label{14}\exists u_0\left.
\begin{array}{rl}
\bigwedge_{i=1}^d u_{0}\cdot r_{i}^{\sigma^{n_{0}-n_{i}}}\equiv_{\tau^{n_{0}-n_{i}}(\delta_{i})} b_{i}\cdot t^{n_{0}-n_{i}}\\
\end{array}\right.
\end{equation}
Suppose we can solve that system. Then by Lemma \ref{actiont}, there exists $u$ such that $u\cdot t^{n_{0}}\equiv_{\delta_{0}}u_{0}$, where we can choose $\delta_{0}\geq \max\{\tau^{n_{0}-n_{i}}(\delta_{i}):\; 1\leq i\leq d\}$.
\par Now we order the set $\{\Upsilon(r_{i}^{\sigma^{n_{0}-n_{i}}},\tau^{n_{0}-n_{i}}(\delta_{i})): 1\leq i\leq d\}$ and we replace
one of the congruences by an equation, the one corresponding to the maximum index. So we are in the case of a separable system treated before.
\par 
\par {\bf Case B':} supppose that for all $1\le i\le d, n_i=0$. We have the system 
\begin{equation}
\label{15}\exists u\left.
\begin{array}{rl}
\bigwedge_{i=1}^d u\cdot r_{i}\equiv_{\delta_{i}} b_{i}\\
\end{array}\right.
\end{equation}
We order the set $\{\Upsilon(r_{i},\delta_{i}): 1\leq i\leq d\}$ and we replace
one of the congruences by an equation, the one corresponding to the maximum index. So we are in the case of a separable system treated before. \qed
\kor\label{eq} In $T_{V}^{sep}$, any $\L_{V}$-formula is equivalent to a quantifier-free $\L_{V}$-formula up to index sentences.
\ekor
\par Recall that index sentences in particular tell us the sizes of the annihilators (of the separable) polynomials and the index of the subgroups $M_{\delta_{1}}.t^{n_{1}}/M_{\delta_{2}}.t^{n_{2}}$, with $\delta_{1}<\delta_{2}$, $n_{1}, n_{2}\in \N$.
Also, the image of $M$ by a $L_{A}$-term $u(x)$ with one free variable $x$ is equal to $M.t^n$, for some $n\in \N$ and we can determine the $n$ from the term $u(x)$, but since our language $\L_{A}$ contains $\lambda$ functions, we need to consider terms in several variables. 
\par In the next section, we will show that if we add a list of axioms specifying the torsion to $T_{V}^{sep}$, then the torsion submodule is determined up to isomorphism (Corollary \ref{iso}). Then in the last section, we will consider the class of torsion-free models of $T_{V}^{sep}$ and we will show that any two elements are elementary equivalent (Corollary \ref{complete}).
\section{Torsion}
\par Let $M\models T_{V}^{sep}$; denote by $M_{tor}$ the submodule of $M$ consisting of torsion elements. Note that $M_{tor}$ is a $\L_{A}$-substructure of $M$ (\cite[Proposition 3.5]{DDP1}. Moreover by our quantifier elimination result (Proposition \ref{EQsep}) $M_{tor}$ is a pure submodule of $M$. So taking an ultrapower $M^*$ of $M$ which is $(\vert A\vert+\aleph_{0})^+$-saturated, the corresponding ultrapower of $M_{tor}$ is a direct summand of $M^*$, namely $M^*=(M_{tor})^*\oplus M_{0}$, where $M_{0}$ is a torsion-free $A$-module, an $\L_{A}$-substructure and a model of $T_{V}^{sep}$.
\par In this section we will show that if we add to the theory $T_{V}^{sep}$ a list of axioms specifying the torsion for each separable polynomial, then the submodule consisting of the torsion elements is unique up to isomorphism as an $\L$-substructure in any model of that extended theory.
\par We will show on one hand that we can determine all the valuations taken by the elements in the annihilator of a separable polynomial belonging of $\I$ and on the other hand that given a non-zero element $n$ of valuation $\delta$ and a separable polynomial  $q(t)\in \I$, we can determine all the valuations taken by the elements $m$ such that $m.q(t)=n$.
\par We will use Corollary \ref{div} together with the factorization of such polynomials $q(t)$ into linear factors of the form $t b-1$, $t-a$, $c$ with $a, b, c\in \O_{K}$, $v(c)=0, v(a)\geq 0, v(b)>0$ (see Proposition \ref{factorization}). 
\nota Let $q(t)\in A$ and let $M$ be an $A$-module, then denote $ann(q(t)):=\{m\in M: m\cdot q(t)=0\}$.
\enota
\lmm\label{unique}
Let $\gamma\in \Gamma$, then there exists at most one $\delta\in \Delta$ such that $\tau(\delta)=\delta+\gamma.$ Moreover if there exists $\delta_{0}$ such that $\tau(\delta_{0})=\delta_{0}+\gamma,$ then 
we have for $\delta>\delta_{0}$ that $\tau(\delta)>\delta+\gamma$ and for $\delta<\delta_{0}$ that $\tau(\delta)<\delta+\gamma$.
\elmm
\pr Let $\delta,\;\delta'\in\Delta$ and suppose that $\tau(\delta)=\delta+\gamma$ and $\tau(\delta')=\delta'+\gamma$ with $\gamma\in \Gamma$. By transitivity of the action of $\Gamma$ on $\Delta$, there exist $\gamma_{1},\;\gamma_{2}\in \Gamma$ such that $\delta=0_{\Delta}+\gamma_{1}$ and $\delta'=0_{\Delta}+\gamma_{2}$. By compatibility of the action of $\tau$ and $\sigma_{v}$, we get that $\tau(\delta)=0_{\Delta}+\sigma_{v}(\gamma_{1})$ and $\tau(\delta')=0_{\Delta}+\sigma_{v}(\gamma_{2})$.
\par We get $\tau(\delta')=\delta'+\gamma=0_{\Delta}+\gamma_{2}+\gamma=0_{\Delta}+\gamma_{1}+(\gamma_{2}-\gamma_{1})+\gamma=\delta+\gamma+(\gamma_{2}-\gamma_{1}).$
So, $\tau(\delta')=\tau(\delta)+(\gamma_{2}-\gamma_{1})$. Therefore,
$0_{\Delta}+\sigma_{v}(\gamma_{2})=0_{\Delta}+\sigma_{v}(\gamma_{1})+(\gamma_{2}-\gamma_{1}).$ Namely, we have $0_{\Delta}+\sigma_{v}(\gamma_{2}-\gamma_{1})=0_{\Delta}+(\gamma_{2}-\gamma_{1}).$ So, $\gamma_{1}=\gamma_{2}$ and $\delta=\delta'$.
\smallskip
\par Now let us show the second part. 
Let $\gamma_{0},\;\gamma_{1}\in \Gamma$ be such that $0_{\Delta}+\gamma_{0}=\delta_{0}$ and $0_{\Delta}+\gamma_{1}=\delta$. Suppose $\delta>\delta_{0}$, then we have that $\gamma_{1}>\gamma_{0}$. By assumption on $\tau$, we have that $\tau(0_{\Delta}+\gamma_{1}-\gamma_{0})>0_{\Delta}+\gamma_{1}-\gamma_{0}$. So, $0_{\Delta}+\sigma_{v}(\gamma_{1}-\gamma_{0})>0_{\Delta}+\gamma_{1}-\gamma_{0}$, equivalently $0_{\Delta}+\sigma_{v}(\gamma_{1})-\sigma_{v}(\gamma_{0})>0_{\Delta}+\gamma_{1}-\gamma_{0}$. So,
$0_{\Delta}+\sigma_{v}(\gamma_{1})>0_{\Delta}+\sigma_{v}(\gamma_{0})+\gamma_{1}-\gamma_{0}=\tau(\delta_{0})+\gamma_{1}-\gamma_{0}=\delta_{0}+\gamma+\gamma_{1}-\gamma_{0}=0_{\Delta}+\gamma+\gamma_{1}=\delta+\gamma$. Namely $\tau(\delta)>\delta+\gamma$.
\par The proof that $\tau(\delta)<\delta+\gamma$ for $\delta<\delta_{0}$ is similar.
\qed
\medskip
\nota Let $\gamma\in \Gamma$ and suppose $\delta\in \Delta$ is such that $\tau(\delta)=\delta+\gamma$, then we will denote $\delta$ by $(\tau-1)^{-1}(\gamma)$. In particular,
$\tau((\tau-1)^{-1}(\gamma))=(\tau-1)^{-1}(\gamma)+\gamma$.

\enota
\lmm\label{linear} Let $M\models T_{V}^{sep}$ and let $r(t)\in \I$ of degree $1$.
\par $(1)$ When $m\in ann(r(t))$, then $w(m)$ takes a unique value which can be expressed in terms of the values of the coefficients of $r(t)$.
\par $(2)$ Let $n\in M-\{0\}$, then there exists $m\in M$ such that $n=m.r(t)$ and $w(m)$ can take at most two values which can be expressed in terms of $w(n)$ and the values of the coefficients of $r(t)$.
\elmm
\pr We can restrict ourselves to consider $r(t)$ of the form $(t-a)$, or $(t b-1)$ with $v(a)\geq 0$ and $v(b)>0$. 
\par (1) Suppose that $m\cdot (t-a)=0$ with $m\neq 0$, then $m\cdot t=m\cdot a$ and so $\tau(w(m))=w(m)+v(a)$. By Lemma \ref{unique}, $w(m)$ is uniquely determined and we will use the above notation:
$(\tau-1)^{-1}(v(a)).$
\par Suppose now that $m\cdot (t b-1)=0$ and $m\neq 0$, then $m\cdot t b=m$ and so $\tau(w(m))+v(b)=w(m)$. We denote $w(m)$ by $(\tau-1)^{-1}(-v(b))$.
\par So in both cases, if there is such a non zero $m$, $w(m)$ can only take one value.
\medskip
\par (2) Now let $n\in M$ with $w(n)=\delta\in \Delta$. By Corollary \ref{div}, there exists $m_{0}$ such that $m_{0}\cdot r(t)=n$ with $w(m_{0})=\Upsilon(r(t),\delta)$ (and any other element $m$ with $m\cdot r(t)=n$ differs from $m_{0}$ by an element of the annihilator of $r(t)$. Let us calculate explicitly $\Upsilon(r(t),\delta)$ in each case. Let $\Upsilon:=\Upsilon(r(t),\delta)$.
\smallskip
\par First, let us consider the case $r(t)=(t b-1)$.
\cl $\;$
\par $(i)$ If $\delta\geq (\tau-1)^{-1}(-v(b))$, then $\Upsilon=\delta$.
\par $(ii)$ If $\delta<(\tau-1)^{-1}(-v(b))$, then $\Upsilon=\tau^{-1}(\delta-v(b))$.
\ecl
\prcl By Lemma \ref{unique}, there is at most one $\rho\in \Delta$ such that $\rho+(-v(b))=\tau(\rho)$ and we have denoted such a $\rho$ by $(\tau-1)^{-1}(-v(b))$. 
\par Moreover if $\rho'<(\tau-1)^{-1}(-v(b))$, then $\rho'>\tau(\rho')+v(b)$ and if $\rho'>(\tau-1)^{-1}(-v(b))$, then $\rho'<\tau(\rho')+v(b)$.
\par We compare $\delta$, respectively $\Upsilon$ to $\tau(\delta)+v(b)$ (equivalently to $(\tau-1)^{-1}(-v(b))$), respectively to $\tau(\Upsilon)+v(b)$.
\par Moreover, by definition $\Upsilon$ is such that $\delta=\min\{\Upsilon,\tau(\Upsilon)+v(b)\}$. 
\par So, if $\Upsilon<\tau(\Upsilon)+v(b)$, then $\delta=\Upsilon$ and this corresponds to the case $\delta<\tau(\delta)+v(b)$ or equivalently to $\delta>(\tau-1)^{-1}(-v(b))$.
\par If $\tau(\Upsilon)+v(b)<\Upsilon$, then $\delta=\tau(\Upsilon)+v(b)$ and so $\delta<\Upsilon$ and so $\tau(\delta)+v(b)<\tau(\Upsilon)+v(b)=\delta$, or equivalently $\delta<(\tau-1)^{-1}(-v(b))$.
\par If $\tau(\Upsilon)+v(b)=\Upsilon$, then $\delta=\Upsilon$ and so $\delta=(\tau-1)^{-1}(-v(b))$.

\qed
\medskip
\par Now suppose there exists $m\neq m_{0}$ such that $m\cdot r(t)=n$, equivalently assume  we have $m_{1}\in ann(r(t))-\{0\}$. Then,
\medskip
\par $(i)$ If $\delta >(\tau-1)^{-1}(-v(b))$, then $w(m_{0}+m_{1})=(\tau-1)^{-1}(-v(b)).$ So in this case we have two possible values for $w(m)$ with $m\cdot (tb-1)=n$. In fact we have one element $m_{0}$ with $w(m_{0})=\delta$ (and $m_{0}\cdot r(t)=n$) and all the other elements $m$ have value $(\tau-1)^{-1}(-v(b))$.
\par $(ii)$ If $\delta <(\tau-1)^{-1}(-v(b))$, then $\tau^{-1}(\delta-v(b))< (\tau-1)^{-1}(-v(b))$. (Indeed, $\rho:=(\tau-1)^{-1}(-v(b))$ is defined by: $\tau(\rho)+v(b)=\rho$. So we have to show that $\delta-v(b)<\tau(\rho)$, equivalently that $\delta<\tau(\rho)+v(b)=\rho.$)
\par So since $\Upsilon=\tau^{-1}(\delta-v(b))$, we have $w(m_{0}+m_{1})=\tau^{-1}(\delta-v(b))$ and in this case we have one possible value for $w(m)$ with $m.r(t)=n$.
\par $(iii)$ If $\delta =(\tau-1)^{-1}(-v(b))$, then $w(m_{0}+m_{1})\geq (\tau-1)^{-1}(-v(b))$. Let us show that we have equality by way of contradiction.
\par Suppose that $w(m_{0}+m_{1})>(\tau-1)^{-1}(-v(b))$. By Lemma \ref{unique}, this implies that
$w(m_{0}+m_{1})<\tau(w(m_{0}+m_{1}))+v(b)$, then $w(n)=\delta=w(m_{0}+m_{1})$
, a contradiction.
\par So, we get that 
$w(m_{0}+m_{1})= (\tau-1)^{-1}(-v(b))=\Upsilon=w(m_{0})$ and again in this case we have only one possible value for $w(m)$ with $m.r(t)=n$.
\medskip
\par Second, let us consider the case $r(t)=(t-a)$. 
\cl $\;$
\par $(i)$ If $\delta\geq (\tau-1)^{-1}(v(a))$, then $\Upsilon=\delta-v(a)$.
\par $(ii)$ If $\delta<(\tau-1)^{-1}(v(a))$, then $\Upsilon=\tau^{-1}(\delta)$.
\ecl
\prcl The proof is similar to the proof of the previous claim. But now, $\Upsilon$ is such that $\delta=\min\{\tau(\Upsilon),\Upsilon+v(a)\}$. 
As before, let $m_{0}$ such that $n=m_{0}\cdot t-m\cdot a$ and $w(m_{0})=\Upsilon$.
So, we have that $\delta \geq \min\{\tau(\Upsilon),\Upsilon+v(a)\}$.
We compare both $\Upsilon$ and $\delta$ to $\tau(\Upsilon)+v(a)$, respectively to $\tau(\delta)+v(a)$ and therefore also to $(\tau-1)^{-1}(v(a))$.\qed
\medskip
\par Again any other solution $m$ of $m\cdot (t+a)=n$ differs from $m_{0}$ by a non zero element $m_{1}$ of the annihilator of $r(t)$. Let us evaluate $w(m_{0}+m_{1})$.
\par $(ia)$ If $(\tau-1)^{-1}(v(a))<\delta<(\tau-1)^{-1}(v(a))+v(a)$. 
\par By the Claim, $w(m_{0})=\delta-v(a)$, and so we have $w(m_{0}+m_{1})=\delta-v(a)$. 
\par $(ib)$ If $\delta >(\tau-1)^{-1}(v(a))+v(a)$, then 
$w(m_{0}+m_{1})=(\tau-1)^{-1}(v(a))=w(m_{1}).$ 
\par $(ic)$ If $\delta =(\tau-1)^{-1}(v(a))+v(a)$, then $w(m_{0}+m_{1})\geq(\tau-1)^{-1}(v(a)).$ Suppose that $w(m_{0}+m_{1})>(\tau-1)^{-1}(v(a))$, so 
$w(m_{0}+m_{1})>\delta-v(a)$. Since $(m_{0}+m_{1})\cdot (t+a)=n$, $w((m_{0}+m_{1})\cdot t)=\delta$. On the other hand, by Lemma \ref{unique},
$\tau(w(m_{0}+m_{1}))>w(m_{0}+m_{1})+v(a)>\delta$, a contradiction. 
So, we also get in this case that $w(m_{0}+m_{1})=(\tau-1)^{-1}(v(a))=w(m_{1}).$
\par So in case $(ia)$, we have two possible values depending on whether there is a non zero element in the annihilator of $r(t)$.
\par $(ii)$ If $\delta <(\tau-1)^{-1}(v(a))$, then by the claim, $w(m_{0})=\tau^{-1}(\delta)$; compare $\tau(w(m_{0}))$ to $\tau(w(m_{1}))$(=$\tau((\tau-1)^{-1}(v(a)))$).
\par We have $\tau((\tau-1)^{-1}(v(a)))=(\tau-1)^{-1}(v(a))+v(a)>\delta$ (see Notation above) and so $w(m_{0})<w(m_{1})$ ($\tau$ respects $<$ on $\Delta$) and so $w(m_{0}+m_{1})=w(m_{0})=\tau^{-1}(\delta)$.
So we have only one possible value.
\par $(iii)$ If $\delta =(\tau-1)^{-1}(v(a))$, then $w(m_{0})=\delta-v(a)=(\tau-1)^{-1}(v(a))-v(a)<(\tau-1)^{-1}(v(a))=w(m_{1})$. Then $w(m_{0}+m_{1})=w(m_{0})=\delta-v(a)$.
So, again in this case we have only one possible value.
\smallskip
\par Note that in each case $w(m)$ can be expressed in terms of $w(n)$ and the values of the coefficients of $r(t)$. \qed
\prop \label{ann} Let $M\models T_{V}^{sep}$, let $m\in M$ and let $q(t)\in \I$ of degree $d$. Then there is a finite subset $F_{q(t)}\subset \Delta$ of cardinality at most $2^{d-1}$ such that if $m\in ann(q(t))-\{0\}$, then $w(m)\in F_{q(t)}$. 
(N.B. The elements of $F_{q(t)}$ whose values are taken by elements of $ann(q(t))$ only depend on the values of the coefficients of the factors of degree $1$ of $q(t)$ and on which are the non-trivial annihilators in $M$.)
\eprop
\pr We proceed by induction on $d$. For polynomials of degree $1$, this is the content of Lemma \ref{linear}. Let us assume $d\geq 2$.
By Proposition \ref{factorization}, $q(t)=r(t) q_{1}(t)$, where $r(t),\;q_{1}(t)\in \I$ and $r(t)$ has degree $1$. Now $m\cdot q(t)=0$ is equivalent to $m\cdot r(t)=0$ or $m\cdot r(t)\in ann(q_{1}(t))-\{0\}$. Since degree of $q_{1}(t)$ is strictly less than $d$, we can apply the induction hypothesis and so we get at most $2^{d-2}$ possible values for the elements in $ann(q_{1}(t))-\{0\}$. By Corollary \ref{div}, for each $n\in ann(q_{1}(t))-\{0\}$, there is an element $m$ such that $m\cdot r(t)=n$ and by Lemma \ref{linear}, for each of the values $w(n)$, we get at most $2$ values for $w(m)$.
\qed
\prop \label{val-div} Let $M\models T_{V}^{sep}$. Let $q(t)\in \I$ of degree $d$ and let $n\in M-\{0\}$.
Then one can determine a finite set $G_{q(t)}\subset \Delta$ of cardinality at most $2^d$ such that $w(m)\in G_{q(t)}$ if and only if $m\cdot q(t)=n$, $m\in M$. 
Moreover, $G_{q(t)}$ only depends on $w(n), q(t)$ and on which are the non-trivial annihilators in $M$.
\eprop
\pr  We proceed by induction on $d\geq 1$.
For polynomials of degree $1$, this is the content of Lemma \ref{linear}. So let us assume $d\geq 2$.
By Proposition \ref{factorization}, $q(t)=q_{2}(t) r(t)$, where $r(t)$ has degree $1$. Now $m\cdot q(t)=n$ is equivalent to $m'\cdot r(t)=n$ and $m\cdot q_{2}(t)=m'$. By Lemma \ref{linear}, given $w(n)$, we know that there is either one or two values for  $w(m')$ with $m'\cdot r(t)=n$ depending on the respective positions of the values of the coefficients of $r(t)$ and $w(n)$, together whether $ann(r(t))$ is non-trivial.
By Corollary \ref{div}, there exists $m$ such that $m\cdot q_{2}(t)=m'$. Then we apply the induction hypothesis to $q_{2}(t)$, so given each of these values for $w(m')$,  the number of values of such $m$ are bounded by $2^{d-1}$ (and we can determine the exact number which depends on the relative position on the chain $\Delta$ of the values of the coefficients and $\delta$ together with which are the non-trivial annihilators).
\qed
\medskip
\par Now we extend $T_{V}^{sep}$ by specifying the torsion in our models. Note that in considering $ann(q(t))$, we may always assume that $q(t)\in \I$, also that annihilators are $Fix(\si)$-vector spaces. If $Fix(\si)$ is infinite then if we have two annihilators with one strictly included in the other, then the index is infinite (\cite[Lemma 2.4 (in Corrigendum)]{BP}). So in this case we will add to the theory $T_{V}^{sep}$ a list of axioms specifying which annihilators are non-trivial.
\par From now on let us assume $Fix(\si)$ is finite. For instance, in the case where $K$ is a (valued) field of characteristic $p$ and $\si$ is the Frobenius endomorphism (or a power of it), then $Fix(\si)$ is finite. We will specify the torsion as follows. 
\defn Let $T_{tor}$ be the theory of $A$-modules together with the following scheme of axioms, for each element $q\in \I$ of degree $d$: there exist exactly $d$ elements $x_{1},\cdots, x_{d}$ which are linearly independent over $Fix(\si)$ and such that $x_{i}\cdot q=0$, for all $1\leq i\leq d$.
\edefn
\par We will consider now the theory $T_{V}^{sep}\cup T_{tor}$. A model of that theory is for instance the separable closure of $K$.

\kor \label{iso} Let $M, N$ be two models of $T_{V}^{sep}\cup T_{tor}$ containing, respectively, isomorphic $\L_{V}$-structures $M_{0},\;N_{0}$. Then we may extend this partial isomorphism to a minimal submodel $M_{0}^{sep}$ of $M$ containing $M_{0}.$
\ekor
\pr We follow the proof of \cite[Proposition 5.8]{DDP1} and we apply Proposition \ref{val-div}.
In particular, a main ingredient in \cite[Proposition 5.8]{DDP1} is the following (see \cite[Lemma 5.1]{DDP1}). Let $N$ be an $\L_{A}$-structure and $N_{0}$ a substructure of $N$. Let $u\in N-N_{0}$ and assume that $u\cdot q\in N_{0}$ for some $q\in \I$.
Then there is a unique (up to multiplication by elements of $Fix(\si)$) element $q_{u}$ of $\I$ such that $u\cdot q_{u}\in N_{0}$. 
\qed
\section{Torsion-free models of $T_{V}^{sep}$}
In this section, we show that the theory of torsion-free models of $T_{V}^{sep}$ is complete, and we specify the other completions of $T_{V}^{sep}$ (see Corollary \ref{complete}).
\par Let $M\models T_{V}^{sep}$ and a torsion-free $A$-module. Given $G_{0},\;G_{1}$ two p.p. definable subgroups of $M$, we wish to determinate the index of $[G_{0}:G_{1}]$. 
\par Here we will assume that the residue field $\bar K$ is infinite, which is the case if $K$ is separably closed of characteristic $p$ and finite (non-zero) imperfection degree, and so $\bar K$ is algebraically closed (and so infinite) and that $\si$ acts on $K$ as the Frobenius.
\par First we note that certain p.p. definable subgroups have infinite index. For instance, the index $[M_{0_{\Delta}}:M_{0_{\Delta}^+}]$ is infinite (we assume that the map $w$ is surjective on $\Delta$, so there is an element $a\in M$ with $w(a)=0$ and multiplying $a$ by elements of $\O_{K}-\M_{K}$ with different images in $\bar K$, we get elements of $M_{0_{\Delta}}$ in distinct cosets modulo $M_{0_{\Delta}^+}$ (and so a fortiori modulo any $M_{\delta}$ with $\delta>0$). Since the action of $\Gamma$ is transitive on $\Delta$, we get the indices $[M_{\delta}:M_{\delta^+}]$ are also infinite, for any $\delta\in \Delta$.
Also the index of $M.t$ in $M$ is infinite (this follows from the fact that $K^{\si}$ is infinite) (\cite[Proposition 3.2]{DDP2}) and the index of $M.t^{m+1}$ in $M.t^m$ as well.
\par Now let us consider the general case. 
\prop \label{infinite} Let $M\models T_{V}^{sep}$ and assume that $M$ is a torsion-free $A$-module. Then the index of any two p.p. definable subgroups $G_{1}\subsetneqq G_{0}$ of $M$ is infinite.
\eprop
\pr
By the positive q.e. result (see Proposition \ref{EQsep}), a p.p. definable subgroup of $M$ is defined by a positive quantifier-free formula of the following kind: $\bigwedge_{i} t_{i}(u)=0\;\&\;\bigwedge_{j} t_{j}(u)\equiv_{\delta_{j}}0,$
where $t_{i},\;t_{j}$ are $\L_{A}$-terms and $u$ one free variable.  As noted before a $\L_{A}$-term $t(u)$ is a $L_{A}$-term in $u,\;\lambda_{i}(u)$, $i\in n^{\ell}$, for some $\ell\geq 1$. Using  the decomposition $t q=\sum_{i}\sqrt[\si]{q_{i}} t c_{i}$ (see Notation 2.2), we may assume that each equation contains a separable coefficient. (This process increases the number of equations but decreases the degree of the coefficients, so it will eventually terminates.)
Note that if we have two equations containing a separable coefficient for $\lambda_{i}(u)$, say $q_{0}$, $q_{1}$, then by multiplying by an element of $K^{\times}$, we may assume that $q_{0}\in \I$ (respectively $q_{1}\in \I$). By the Ore property of $A_{0}$, there exist $q_{0}', \;q_{1}'$ such that $q_{0} q_{1}'=q_{1} q_{0}'$ and note that we may choose $q_{0}', q_{1}'\in \I$. So we may just keep one equation with $\lambda_{i}(u)$ and assume that it has a separable coefficient. Continuing like that we get a lower triangular $d\times k$-matrix $D$ with non zero separable coefficients on its diagonal. Recall that such matrix was called {\it lower triangular separable} (l.t.s.) of co-rank $\ell$ in \cite[Definition 9]{DDP1} where $d-\ell$ is the number of non-zero (or equivalently separable) elements on its diagonal. We will call the corresponding system a l.t.s. system.
\par So we may assume that $G_{0}$ (respectively $G_{1}$) is defined by a l.t.s. system of co-rank $\ell\leq d$ (respectively $\ell_{1}\leq \ell$) equations (on the same variables) with in addition congruences conditions. 
\par Set $u_{d-\ell}:=\lambda_{m}(u)$, where $\lambda_{m}(u)$ is the variable multiplied by the $(d-\ell,d-\ell)$-coefficient of $D$.
Let $\delta_{i}$, $i\in m$, be the elements of $\Delta$ occurring in the congruences. 
\par Since $G_{1}\subseteq G_{0}$, we may assume that $G_{1}$ is defined by a l.t.s. system of co-rank $\ell_{1}\leq \ell$ and with the same first $d-\ell$ equations. 
\cl\label{tf} Let $q\in \I$, then in any torsion-free model of $T_{V}^{sep}$, we have the following equivalence: $u\equiv_{\delta} 0\leftrightarrow u\cdot t^n q\equiv_{\Upsilon^{-1}(q,\tau^n(\delta))} 0$.
\ecl
\prcl By definition of $\Upsilon^{-1}$ (see Notation \ref{notation-upsilon}), we get in any valued $A$-module, that $u\equiv_{\delta} 0\rightarrow u\cdot t^n q\equiv_{\Upsilon^{-1}(q,\tau^n(\delta))} 0$.
Now assume that $w(u\cdot t^n q)\geq \Upsilon^{-1}(q,\tau^n(\delta))$, since $q\in \I$, there exists $u'$ with $w(u')\geq \tau^n(\delta)$ ($\Upsilon$ is increasing, see Lemma \ref{angularbis}) such that $u'\cdot q=u\cdot t^n q$. Since $M$ is torsion-free, $u'=u\cdot t^n$ and so $w(u)\geq \delta$. \qed
\par So if we have a l.t.s. system of co-rank $\ell\leq d$, with $d\geq 1$ and if $u_{m}$, $1\leq m\leq d-\ell$, occurs in an equation with a separable coefficient and if it also occurs in a congruence condition, we may replace it in the congruences conditions in terms of the other variables (using the Ore property of $A_{0}$ and the Claim above). So w.l.o.g., we may assume that in the congruences none of the variables $u_{m}$, $1\leq m\leq d-\ell$ occur.

\par First we assume that $\ell_{1}<\ell$, so we have at least one more equation in $G_{1}$ of the form $u_{d-\ell+1}\cdot q_{d-\ell+1}+\cdots=0$ $(\star)$.
Denote by $\delta(G_{0})$ the minimum $\delta\in \Delta$ such that any tuple of elements each in $M_{\delta}$ satisfies the congruence conditions appearing in the definition of $G_{0}$ and denote by $\chi$ the conjunction of these congruences.
\par Take $\overrightarrow{u}=(u_{d-\ell+1},\overrightarrow{u_{0}})$ satisfying $\chi$ and 
verify whether $\overrightarrow{u}$ satisfies equation $(\star)$. If it does, add to $u_{d-\ell+1}$ a non-zero element $u_{\delta,1}\in M_{\delta(G_{0})}$. So the tuple $(u_{n-\ell+1}+u_{\delta,1},\overrightarrow{u_{0}})$ still satisfies $\chi$ (since $(u_{\delta,1},\overrightarrow{0})$ satisfies $\chi$) but no longer equation $(\star)$ ($u_{d-\ell+1}$ was uniquely determined in terms of $\overrightarrow{u_{0}}$ since we are in a torsion-free module). 
Since the system is l.t.s. and $M\models T_{V}^{sep}$, we may find
$(u_{1},\cdots,u_{d-\ell})$ such that $(u_{1},\cdots,u_{d-\ell},u_{d-\ell+1}+u_{\delta,1},\overrightarrow{u_{0}})\in G_{0}$, then choose in $M_{\delta(G_{0})}$ a non zero $u_{\delta,2}\neq u_{\delta,1}$. The tuple $(u_{d-\ell+1}+u_{\delta,2},\overrightarrow{u_{0}})$ still satisfies $\chi$ but no longer equation $(\star)$. Again we may complete this tuple to $(u_{1}',\cdots,u_{d-\ell}',u_{d-\ell+1}+u_{\delta,2},\overrightarrow{u_{0}})\in G_{0},$ and so $(u_{1}'-u_{1},\cdots,u_{d-\ell}'-u_{d-\ell},u_{\delta,2}-u_{\delta,1},\overrightarrow{0})\in G_{0}-G_{1}$, (since $(u_{\delta,2}-u_{\delta,1},\overrightarrow{0})$ does not satisfy equation $(\star)$) and we may continue infinitely often since $M_{\delta(G_{0})}$ is infinite.
\par Then assume that $\ell_{1}=\ell$ and so that the l.t.s. system of equations occurring in the definition of $G_{0}$ and $G_{1}$ is the same.
\par Then using the Claim \ref{tf}, we triangularize further the system of congruences as follows.
\par Suppose we have the following system of two congruences with $r_{1}, r_{2}\in \I$ and $\delta_{1}\geq \delta_{2}\in \Delta$, occurring in the definition of $G_{0}$ and for simplicity re-index the variables by $u_{0}, u_{1},\cdots$.
\begin{equation}
\label{congruence}\left.
\begin{array}{rl}
u_{0}\cdot r_{1}+u_{1}\cdot r_{3}+\cdots \equiv_{\delta_{1}} 0\\
u_{0}\cdot r_{2}+u_{1}\cdot r_{4}+\cdots \equiv_{\delta_{2}} 0
\end{array}\right.
\end{equation}
By the right Ore property of $A_{0}$, there exist $q_{1}, q_{2}\in A_{0}$ (and we may assume that they belong to $\I$) such that $r_{1}q_{2}=r_{2}q_{1}$. 
\par Assume $\Upsilon^{-1}(q_{2},\delta_{1})\leq \Upsilon^{-1}(q_{1},\delta_{2})$ (the other case is similar). Then system (\ref{congruence}) is equivalent to:
\begin{equation}
\label{congruence2}\left.
\begin{array}{rl}
u_{0}\cdot r_{2}+u_{1}\cdot r_{4}+\cdots \equiv_{\delta_{2}} 0\\
u_{1}\cdot (r_{3}q_{2}-r_{4}q_{1})+\cdots \equiv_{\Upsilon^{-1}(q_{2},\delta_{1})} 0
\end{array}\right.
\end{equation}
Given a solution of the system (\ref{congruence2}), we multiply the first equation by $q_{1}$ and get
$u_{0}\cdot r_{2}q_{1}+u_{1}\cdot r_{4}q_{1}+\cdots \equiv_{\Upsilon^{-1}(q_{1},\delta_{2})} 0$ and add the second equation to get 
$u_{0}\cdot r_{1}q_{2}+u_{1}\cdot r_{3}q_{2}+\cdots \equiv_{\Upsilon^{-1}(q_{2},\delta_{1})} 0$. This last equation is equivalent by Claim \ref{tf} to 
$u_{0}\cdot r_{1}+u_{1}\cdot r_{3}+\cdots \equiv_{\delta_{1}} 0$.
\par So, w.l.o.g. we may assume that we have the following system of two congruences with $r_{1}, r_{2}\in \I$, the first one being the last one occurring in the definition of $G_{0}$ and the second one being the first one in $G_{1}$ and for simplicity re-index the variables by $u_{0}, u_{1},\cdots$.
\begin{equation}
\label{congruence3}\left.
\begin{array}{rl}
u_{0}\cdot r_{1}+u_{1}\cdot r_{3}+\cdots \equiv_{\delta_{1}} 0\\
u_{1}\cdot r_{4}+\cdots \equiv_{\delta_{2}} 0
\end{array}\right.
\end{equation}
First choose $(u_{1},\overrightarrow{0})$ such that the second equation does not hold (take $u_{1}\notin M_{\Upsilon^{-1}(r_{4},\delta_{2})}$). Then choose 
$u_{0}$ such that $u_{0}\cdot r_{1}+u_{1}\cdot r_{3}+\cdots = 0$, which is possible whenever $r_{1}$ is separable and in the case $r_{1}$ is of the form $t^m.r_{1}'$ with $r_{1}'\in \I$ separable and $m\geq 1$, we proceed as follows. First choose, $u_{0}'$ such that 
$u_{0}'\cdot r_{1}'+u_{1}\cdot r_{3}+\cdots = 0$ holds and then $u_{0}$ with $w(u_{0}-u_{0}')\geq \delta_{1}$ (see Corollary \ref{actiont}). So in both cases, $(u_{0},u_{1},\overrightarrow{0})$ satisfies (\ref{congruence3}) and also the conjunction $\chi$ of congruences occurring in the definition of $G_{0}$ but doesn't satisfy those occurring in the definition of $G_{1}$.
Then continue to solve the system in order to get an element in $G_{0}$ using the fact we have put the system in triangular form.
Finally choose infinitely many such $u_{1}$ not congruent modulo $M_{\Upsilon^{-1}(r_{4},\delta_{2})}$.\qed
\kor \label{complete} The theory of torsion-free (as $A$-modules) models of $T_{V}^{sep}$ is complete.
\par Suppose $Fix(\si)$ is finite. Then the theory $T_{V}^{sep}\cup T_{tor}$ is complete.
\par Suppose that $Fix(\si)$ is infinite. Then any extension of $T_{V}^{sep}$ containing the list of axioms telling for each separable $q\in \I$ whether the annihilator of that element is trivial or not, is complete.
\ekor
\pr For the first statement, apply Proposition \ref{eq} and Proposition \ref{infinite}. For the second statement, apply Proposition \ref{eq},  Corollary \ref{iso} and Proposition \ref{infinite}.  And for the third one apply
Proposition \ref{eq} and Proposition \ref{infinite}. 
\qed



\begin{thebibliography}{9}
\bibitem{Az} Azgin S., Valued fields with contractive automorphisms and Kaplansky fields, J. Algebra 324 (2010) 2727-2785.
\bibitem{BP} B\'elair L., Point F., Quantifier elimination in valued Ore modules, J. Symb. Logic 75 
 (2010) 1007-1034. Corrigendum : J. Symb. Logic 77 (2012) 727-728.
\bibitem{C} Cohn P.M., Skew fields, {Encyclopedia of mathematics and its applications} (G.-C. Rota, Editor), vol. 57, Cambridge University Press, 1995.
\bibitem{Co} Conrad P.F., Embedding Theorems for Abelian Groups with Valuations, American Journal of Mathematics 75 (1953), no. 1, 1-29.
\bibitem{DDP1} P. Dellunde, F. Delon, F. Point, The theory of modules
of separably closed fields 1, Journal of Symbolic Logic 67  (2002),  
no. 3, 997-1015.
\bibitem{DDP2} Dellunde, P., Delon F., Point F., The theory of modules
 of separably closed fields-2, Ann. Pure Appl. Logic  129  (2004),  no. 1-3, 181-210.
\bibitem{Del} Delon F., Quelques propri\'et\'es des corps valu\'es en th\'eorie des mod\`eles, Th\`ese d'\'etat Paris 7, 1982.
\bibitem{D} Delon F., Id\'eaux et types sur les corps s\'eparablement clos, M\'em. Soc. Math. France (N.S.) No. 33 (1988), 76 pages.
\bibitem{DS} Denef, J., Schoutens, H., On the decidability of the existential theory of ${\Bbb F}_p[\![t]\!]$, Valuation theory and its applications, Vol. II (Saskatoon, SK, 1999), 43-60, Fields Inst. Commun., 33, Amer. Math. Soc., Providence, RI, 2003. 
\bibitem{Dries} van den Dries L., Quantifier elimination for linear formulas over ordered and valued fields, in: {\em Proceedings of the Model Theory Meeting (Univ. Brussels, Brussels/Univ. Mons, Mons, 1980)},  Bull. Soc. Math. Belg. SŽr. B  33  (1981), no. 1, pp. 19-31.
\bibitem{F} Fleischer I., Maximality and ultracompleteness in valued normed modules, Proc. Amer. Math. Soc., volume 9, number 1, 1958, 151-157.
\bibitem{H} Haran D., Quantifier elimination in separably closed fields of finite imperfection degree, Journal of Symbolic Logic, vol. 53, no.2, 1988, 463-469.
\bibitem{Ho} Hong J., Immediate expansions by valuations of fields, PhD thesis McMaster, august  2013.
\bibitem{L} Lazard M., Groupes $p$-adiques analytiques, Inst. Hautes Etudes Sci. Publ. Math. 26, 1965, 389-603.
\bibitem{On} Onay  G., Modules valu\'es, PhD Thesis Paris 7, december 2011.
\bibitem{Ore} Ore O., On a special class of polynomials, Trans. Amer.
Math. Soc. 35, pp. 559-584, 1933.
\bibitem{Ore2}  Ore O., Theory of non-commutative polynomials, Ann.
Math. 34, pp. 480-508, 1933.
\bibitem{Pa} Pal K., Multiplicative valued difference fields, The Journal of Symbolic Logic 77, number 2, June 2012, 545-579.
\bibitem{P} Point F., Asymptotic theory of modules of separably closed fields, J. Symb. Logic 70 (2005) 573-592.
\bibitem{Pr} Prest M., The model theory of modules, L.M.Soc. Lecture Notes
Series 130, Cambridge University Press, 1988.
\bibitem{R} Rohwer T.,  Valued difference fields as modules over twisted polynomial rings, PhD thesis, University of Illinois at Urbana-Champaign, 2003.
\bibitem{Sr} Srour G., The independence relation in separably closed fields, Journal of Symbolic Logic, vol. 51, 1986, 715-725.

\end{thebibliography}
\end{document}